\documentclass[11pt,a4paper]{article}
\usepackage{color,amsmath,amssymb, lscape,graphicx,ifthen,pdfpages}
\DeclareSymbolFont{boldsymbols}{OMS}{cmsy}{b}{n} 
\DeclareSymbolFontAlphabet{\mathbfcal}{boldsymbols} 

\usepackage[colorlinks=true,
linkcolor=webgreen,
filecolor=webbrown,
citecolor=webgreen]{hyperref}
\definecolor{webgreen}{rgb}{0,.5,0}\definecolor{webbrown}{rgb}{.6,0,0}

\newcommand{\seqnum}[1]{\href{http://oeis.org/#1}{\underline{#1}}}

\def\dstyle#1{$\displaystyle #1 $}

\def\pn{\par\noindent}

\def\psn{\par\smallskip\noindent}
\def\pbn{\par\bigskip\noindent}
\def\Beq{\begin{equation}}
\def\Eeq{\end{equation}}
\def\Beqarray{\begin{eqnarray}}
\def\Eeqarray{\end{eqnarray}}

\def\sspkl{\, < \,}
\def\sspgr{\, > \,}
\def\sspeq{\, =\,}
\def\speq{\ =\ }
\def\sspdef{\, :=\,}

\def\sspfed{\, =:\,}
\def\sspin{\, \in \,}

\def\sspp{\, +\ }
\def\sspm{\, -\ }

\def\sspto{\,\to\,}

\def\sspneq{\, \neq \,}
\def\sspequiv{\,\equiv\,}

\def\ie{{\it i.e.},\, }
\def\eg{{\it e.g.},\, }

\def\viz{{\it viz}\, }

\def\sspfed{\, =:\, }

\def\abs#1{\vert {\,#1\,} \vert}
\def\sspapprox{\,\approx\,}
\def\via{{\it via}\,}

\definecolor{GreenYellow}{cmyk}{0.15,0,0.69,0}
\definecolor{Yellow}{cmyk}{0,0,1,0}
\definecolor{Goldenrod}{cmyk}{0,0.10,0.84,0}
\definecolor{Dandelion}{cmyk}{0,0.29,0.84,0}
\definecolor{Apricot}{cmyk}{0,0.32,0.52,0}
\definecolor{Peach}{cmyk}{0,0.50,0.70,0}
\definecolor{Melon}{cmyk}{0,0.46,0.50,0}
\definecolor{YellowOrange}{cmyk}{0,0.42,1,0}
\definecolor{Orange}{cmyk}{0,0.61,0.87,0}
\definecolor{BurntOrange}{cmyk}{0,0.51,1,0}
\definecolor{Bittersweet}{cmyk}{0,0.75,1,0.24}
\definecolor{RedOrange}{cmyk}{0,0.77,0.87,0}
\definecolor{Mahogany}{cmyk}{0,0.85,0.87,0.35}
\definecolor{Maroon}{cmyk}{0,0.87,0.68,0.32}
\definecolor{BrickRed}{cmyk}{0,0.89,0.94,0.28}
\definecolor{Red}{cmyk}{0,1,1,0}
\definecolor{OrangeRed}{cmyk}{0,1,0.50,0}
\definecolor{RubineRed}{cmyk}{0,1,0.13,0}
\definecolor{WildStrawberry}{cmyk}{0,0.96,0.39,0}
\definecolor{Salmon}{cmyk}{0,0.53,0.38,0}
\definecolor{CarnationPink}{cmyk}{0,0.63,0,0}
\definecolor{Magenta}{cmyk}{0,1,0,0}
\definecolor{VioletRed}{cmyk}{0,0.81,0,0}
\definecolor{Rhodamine}{cmyk}{0,0.82,0,0}
\definecolor{Mulberry}{cmyk}{0.34,0.90,0,0.02}
\definecolor{RedViolet}{cmyk}{0.07,0.90,0,0.34}
\definecolor{Fuchsia}{cmyk}{0.47,0.91,0,0.08}
\definecolor{Lavender}{cmyk}{0,0.48,0,0}
\definecolor{Thistle}{cmyk}{0.12,0.59,0,0}
\definecolor{Orchid}{cmyk}{0.32,0.64,0,0}
\definecolor{DarkOrchid}{cmyk}{0.40,0.80,0.20,0}
\definecolor{Purple}{cmyk}{0.45,0.86,0,0}
\definecolor{Plum}{cmyk}{0.50,1,0,0}
\definecolor{Violet}{cmyk}{0.79,0.88,0,0}
\definecolor{RoyalPurple}{cmyk}{0.75,0.90,0,0}
\definecolor{BlueViolet}{cmyk}{0.86,0.91,0,0.04}
\definecolor{Periwinkle}{cmyk}{0.57,0.55,0,0}
\definecolor{CadetBlue}{cmyk}{0.62,0.57,0.23,0}
\definecolor{CornflowerBlue}{cmyk}{0.65,0.13,0,0}
\definecolor{MidnightBlue}{cmyk}{0.98,0.13,0,0.43}
\definecolor{NavyBlue}{cmyk}{0.94,0.54,0,0}
\definecolor{RoyalBlue}{cmyk}{1,0.50,0,0}
\definecolor{Blue}{cmyk}{1,1,0,0}
\definecolor{Cerulean}{cmyk}{0.94,0.11,0,0}
\definecolor{Cyan}{cmyk}{1,0,0,0}
\definecolor{ProcessBlue}{cmyk}{0.96,0,0,0}
\definecolor{SkyBlue}{cmyk}{0.62,0,0.12,0}
\definecolor{Turquoise}{cmyk}{0.85,0,0.20,0}
\definecolor{TealBlue}{cmyk}{0.86,0,0.34,0.02}
\definecolor{Aquamarine}{cmyk}{0.82,0,0.30,0}
\definecolor{BlueGreen}{cmyk}{0.85,0,0.33,0}
\definecolor{Emerald}{cmyk}{1,0,0.50,0}
\definecolor{JungleGreen}{cmyk}{0.99,0,0.52,0}
\definecolor{SeaGreen}{cmyk}{0.69,0,0.50,0}
\definecolor{Green}{cmyk}{1,0,1,0}
\definecolor{ForestGreen}{cmyk}{0.91,0,0.88,0.12}
\definecolor{PineGreen}{cmyk}{0.92,0,0.59,0.25}
\definecolor{LimeGreen}{cmyk}{0.50,0,1,0}
\definecolor{YellowGreen}{cmyk}{0.44,0,0.74,0}
\definecolor{SpringGreen}{cmyk}{0.26,0,0.76,0}
\definecolor{OliveGreen}{cmyk}{0.64,0,0.95,0.40}
\definecolor{RawSienna}{cmyk}{0,0.72,1,0.45}
\definecolor{Sepia}{cmyk}{0,0.83,1,0.70}
\definecolor{Brown}{cmyk}{0,0.81,1,0.60}
\definecolor{Tan}{cmyk}{0.14,0.42,0.56,0}
\definecolor{Gray}{cmyk}{0,0,0,0.50}
\definecolor{Black}{cmyk}{0,0,0,1}
\definecolor{White}{cmyk}{0,0,0,0}

\normalbaselines
\begin{document}
\bibliographystyle{unsrt}
\rightline{Karlsruhe} \par\smallskip\noindent
\rightline{September 02  2014}
\vbox {\vspace{6mm}}
\begin{center}
{\Large {\bf Notes on Some Geometric and Algebraic Problems Solved by Origami }}\\ [9mm]
Wolfdieter L a n g \footnote{ 
\href{mailto:wolfdieter.lang@partner.kit.edu}{\tt wolfdieter.lang@partner.kit.edu},\quad 
\url{http://www-itp.particle.uni-karlsruhe.de/~wl}
                                          } \\[3mm]
\end{center}
\vspace{2mm}
\centerline {\bf Abstract}\psn
Details for known solutions of some geometric and algebraic problems with the help of origami are presented: two theorems of Haga, the general cubic equation, especially the heptagon equation, doubling the cube as well as the trisection of  angles $\alpha$,  $\pi\sspm \alpha$ and $\pi\sspp \alpha$.\pbn
{\bf Introductory remarks}\psn
These notes give details on some geometric and algebraic problems related to cubic equations which are solved using origami (Japanese for folding paper). Seven axioms for origami can be found in \cite{Axioms}.\pn
These notes mostly start with a square of some given length, called $R$ in some length unit. Given any (transparent) sheet larger than this square one can fold an $R\times R$ square, provided one can determine the distance $R$ between two points $P$ and $Q$ on some line (crease). This assumes that one has some way to measure $R$, \eg a marked ruler. Then one starts with some crease, call it  $c_1$, and folds perpendicular to this crease through some point, defined as the first corner $A$, another crease called $c_2$. A perpendicular folding with respect to some line $c$ (crease) and a point $P$ (not necessarily on $c$) can be accomplished (sometimes called axiom 4 or IV), but here it is useful to have a transparent sheet in order to see when the two parts of $c$  fit together for this folding through $P$. Then one finds the next corner of the square, called $D$, at a given distance $R$ from $A$ on the crease $c_1$. Next, through $D$ a crease $c_3$ perpendicular to $c_1$ is formed. The next crease $c_4$ is obtained by folding crease $c_1$ onto crease $c_2$ (point $A$ will lie on both creases; guaranteed by axiom 3 or III). This will define the next corner $C$ as the intersection point of $c_4$ with $c_3$. Finally a crease $c_5$  perpendicular to $c_3$ through point $C$ is formed to find the last corner $B$ as the intersection of $c_5$ with $c_2$. Alternatively one can fold crease $c_1$ onto crease $c_3$, with $D$ on both creases, to find $B$ as the intersection with crease $c_2$.  This completes then the square $A,\, B,\, C,\, D$ oriented in the positive sense (on one of the transparent paper's sides). In the following we will use the notation $\overline{B,C}$ to denote the straight line connecting points $B$ and $C$, as well as the length of this line segment. The latter should be denoted by $\abs{\overline{B,C}}$, but it should be clear what is meant in each case.\pbn
{\bf Problem 1:\ Haga's second Theorem}\psn
 In the book of {\sl Bellos} \cite{Bellos} one finds on p. 115 an origami leading to the ``second theorem'' of {\sl Kazuo Haga}. For this one folds two neighboring corners of a square sheet of paper (length of the side $R$ in some unit), say $B$ and $C$, in turn on some point $B^{\prime}\sspeq C^{\prime}$ of the opposite side (bordered by the corners $A$ and $D$). This is done by two intersecting creases (called $f$ and $g$ in {\sl Figure 3}). The intersection point, called  $S$ in {\sl Figure 3}, will always lie on the crease which arises if one folds $C$ onto $B$ (which gives one of the medians of the square). This happens independently of the position of the point on which the two corners have been folded. In addition, the three distances between the intersection point $S$ and the chosen point  $C^{\prime}$ and the two corners $B$  and $C$ coincide. \pn
In order to analyse this consider first {\it Figure 1} where $C$ is folded onto $C^{\prime}$ having distance $x\,R$ from the left upper corner $A$. \psn

\parbox{16cm}{
{\includegraphics[height=8cm,width=.5\linewidth]{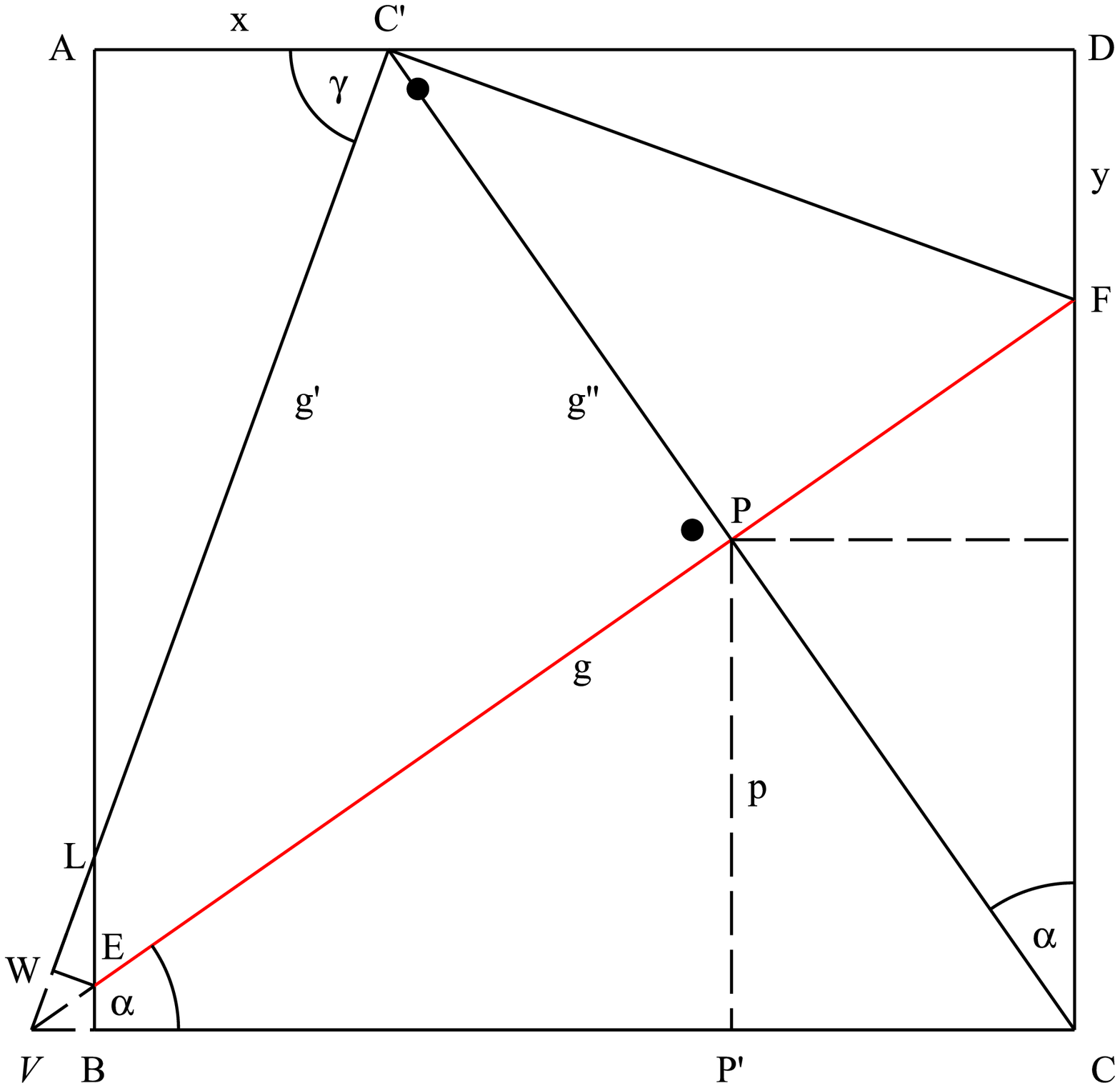}}
{\includegraphics[height=8cm,width=.5\linewidth]{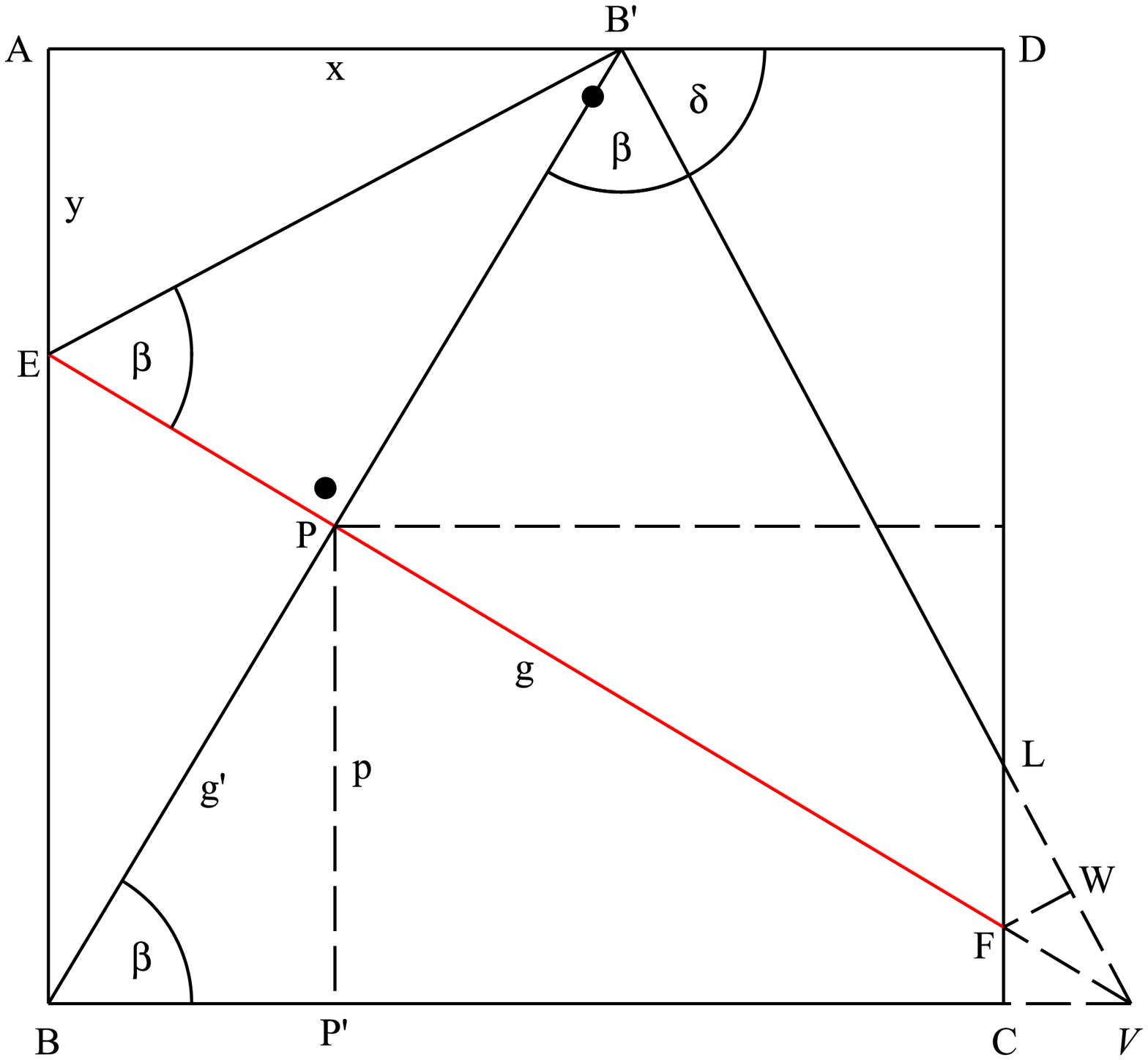}}
}
\psn
\hskip 2.5cm  Figure 1:\  Folding $C$ onto $C^\prime$ \hskip 3cm Figure 2:\  Folding $B$ onto $B^\prime$\pbn
{\bf Figure 1}:\psn
$\overline{A,B}\sspeq  R \sspeq  \overline{A,D},\ \overline{A ,C^{\prime}}\sspeq x\,R,\ \overline{D,F}\sspeq y\,R,\ \overline{V,B}\sspeq \overline{V,W}\sspeq v\,R,\ \overline{W,C^{\prime}}\sspeq \overline{B,C}\sspeq R,\   \overline{B,E}\sspeq \overline{W,E}\sspeq z\,R,\  \overline{P,C}\sspeq \overline{P,C^{\prime}}\sspeq c\,R,\ \overline{P,P^{\prime}}\sspeq p\,R,\  \angle(F,V,C)\sspeq  \angle(C^{\prime},V,F) \sspeq  \angle(C^{\prime},C,D)\sspeq \angle(P,C^{\prime},F) \sspeq \alpha,\  \angle(V,P,C^{\prime})\sspeq \frac{\pi}{2}\sspeq  \angle(V,C^{\prime},F)\ $ (indicated by the two bullets), \ $\angle(V,C^{\prime},A)\sspeq \gamma \sspeq 2\,\alpha$. \pn
$P$ is the intersection point of two perpendicular straight lines, \viz\ $g''$ connecting $C^{\prime}$ and $C$ and $g$ (the crease bringing $C$ to $C^{\prime}$) connecting $V$ and $F$.\pn
The analytic data, depending on $R$ (usually taken as $1$ length unit) and $x$, is :\pn
\dstyle{y\sspeq \frac{x\,(x-2)}{2},\ 2\,c\sspeq \sqrt{1\sspp (1-x)^2},\ \overline{P,F}\sspeq R\,\sqrt{(1-y)^2\sspm c^2},\ \tan\, \alpha\sspeq \frac{\overline{F,C}}{\overline{V,C}} \sspeq \frac{z}{v} \sspeq \frac{\overline{P,F}}{c\,R}\sspeq \frac{\overline{C^{\prime},D}}{R}  \sspeq 1-x,\ \tan(2\,\alpha)\sspeq \frac{2\,(1-x)}{x\,(1-x)},\ \sin\,\alpha\sspeq \frac{1-x}{\sqrt{1\sspp(1-x)^2}}\sspeq \frac{c}{1+v},\ \cos\,\alpha\sspeq \frac{1}{\sqrt{1\sspp(1-x)^2}}\sspeq \frac{p}{c},\  p\sspeq \frac{1}{2},\   v\sspeq \frac{x^2}{2\,(1-x)},\ \overline{F,C}\sspeq \overline{F,C^{\prime}}\sspeq R\,(1\sspm y)\sspeq R\,\frac{2-2\,x-x^2}{2},\  z\sspeq v\,\tan\,\alpha\sspeq \frac{x^2}{2},\ \overline{V,E}\sspeq R\,\sqrt{z^2\sspp v^2}\sspeq R\,\frac{x^2\,\sqrt{1\sspp(1-x)^2}}{2\,(1-x)},\ \overline{W,L}\sspeq R\,z\,\tan(2\,\alpha)\sspeq R\,\frac{x^2\,(1-x)}{\sqrt{x\,(2-x)}},\ \overline{L,A}\sspeq R\,x\,\tan(2\,\alpha)\sspeq R\,\frac{2\,(1-x)}{2-x},\  \overline{E,L}\sspeq\frac{z\,R}{\cos(2\,\alpha)}\sspeq R\,\frac{x\,(1\sspp (1-x)^2)}{2-x},\ \overline{L,C^{\prime}}\sspeq\ R\, \frac{x}{\sin(\frac{\pi}{2} \sspm 2\,\alpha)} \sspeq R\,\frac{1\sspp (1-x)^2}{2-x},\ \overline{L,B}\sspeq\ R\,v\, \tan(2\,\alpha) \sspeq R\,\frac{x}{2-x}, \  \overline{E,P}\sspeq \frac{p\,R}{sin\,\alpha}\sspm \overline{V,E} \sspeq R\,\frac{(1+x)\,\sqrt{1\sspp (1 \sspm x)^2}}{2}},\pn
\dstyle{\overline{B,P^{\prime}}\sspeq R\,(1\sspm c\, \sin\,\alpha)\sspeq R\, \frac{1+x}{2},\ \overline{P,C^{\prime}}\sspeq R\,\frac{1-x}{2}}. \pn
With the origin $O\sspeq B$ in the $(\hat x,\hat y)$-plane (no confusion with the above $x$ and $y$ should arise) the straight lines $g, g'$ and $g''$ are given by:\pn
\dstyle{g:\ \hat y \sspeq (1\sspm x)\,(\hat x \sspp v\,R),\ g':\ \hat y\sspeq \frac{2\,(1-x)}{x\,(2-x)}\,(\hat x\sspp v\,R),\ \ g'':\ \hat y\sspeq -\frac{1}{1-x}\,(\hat x-R)\, .}\psn
{\bf Observation $\bf 1$:} If $x$ varies from $0$ to $R$ then $P$ moves on the middle line \dstyle{\hat y \sspeq \frac{R}{2}} from \dstyle{\hat x\sspeq \frac{R}{2}} to $R$.
\pbn
Next, the analysis is done for the case when the left lower corner $B$ is folded onto $B^{\prime}$ on the side connecting the corners $A$ and $D$, with distance $x\,R$ from $A$ (not necessarily the same $x$ as in {\it Figure 1}). \psn
{\bf Figure 2}:
\psn
$\overline{A,B}\sspeq  R \sspeq  \overline{A,D},\  \overline{A ,B^{\prime}}\sspeq x\,R,\, \overline{A,E}\sspeq y\,R,\  \overline{V,C}\sspeq \overline{V,W}\sspeq v\,R,\ \overline{E,B}\sspeq \overline{E,B^{\prime}}\sspeq (1-y)\,R,\ \overline{W,B^{\prime}}\sspeq \overline{C,B}\sspeq R,\ \overline{C,F}\sspeq \overline{F,W}\sspeq z\,R,\  \overline{P,B}\sspeq \overline{P,B^{\prime}}\sspeq b\,R,\ \overline{P,P^{\prime}}\sspeq p\,R,\  \angle(P,B,V)\sspeq  \angle(P,B^{\prime},V) \sspeq \angle(A,B^{\prime},B) \sspeq \angle(B^{\prime},E,F)\sspeq \angle(P,E,B) \sspeq \angle(P^{\prime},P,V) \sspeq \beta,\  \angle(L,B^{\prime},D)\sspeq \delta \sspeq \pi\sspm 2\,\beta,\ \angle(C,L,V)\sspeq \angle(B^{\prime},L,D)\sspeq  \angle(E,B^{\prime},A)  \sspeq 2\,\beta\sspm \frac{\pi}{2}$ . \pn
 $P$  is the intersection point of the perpendicular straight lines $g'$, with points $B^{\prime}$ and $B$, and $g$ (the crease bringing $B$ to $B^{\prime}$) with points $V$ and $E$.\pn
The analytic data which depends on $R$ (usually taken as $1$ length unit) and $x$ is:\pn
\dstyle{y\sspeq \frac{1-x^2}{2},\ 2\,b\sspeq \sqrt{1\sspp x^2},\ \overline{P,E}\sspeq R\,\sqrt{(1-y)^2\sspm b^2},\ \tan (\frac{\pi}{2}\sspm\beta)\sspeq \frac{1}{\tan\,\beta} \sspeq \frac{\overline{E,B}}{\overline{B,V}} \sspeq \frac{\overline{E,B^{\prime}}}{\overline{B^{\prime},V}} \sspeq  \frac{z}{v} \sspeq \frac{b\,R}{\overline{P,V}} \sspeq  \frac{\overline{E,P}}{b\,R} \sspeq  \frac{\overline{A,B^{\prime}}}{\overline{A,B}} \sspeq x,\ \tan(2\,\beta-\frac{\pi}{2})\sspeq -\frac{1}{\tan(2\,\beta)}\sspeq \frac{y}{x}\sspeq \frac{1-x^2}{2\,x},\   \cos\,\beta\sspeq \frac{x}{\sqrt{1\sspp x^2}}\sspeq  \frac{\overline{B,P^{\prime}} }{b\,R}\sspeq \frac{b}{1+v},\ \overline{B,P^{\prime}}\sspeq R\,\frac{x}{2},\  \ v\sspeq \frac{(1-x)^2}{2\,x},\ \sin\,\beta\sspeq \frac{1}{\sqrt{1\sspp x^2}}\sspeq  \frac{p}{b},\  p\sspeq \frac{1}{2},\ \overline{E,B}\sspeq \overline{E,B^{\prime}}\sspeq R\,(1\sspm y)\sspeq R\,\frac{1+x^2}{2},\  z\sspeq v/\tan\,\beta \sspeq \frac{(1-x)^2}{2},\ \overline{V,F}\sspeq R\,\sqrt{z^2\sspp v^2}\sspeq R\,\frac{(1-x)^2}{2\,x}\, \sqrt{1+x^2},\ \overline{W,L}\sspeq R\,\frac{z}{\tan(2\,\beta\sspm\frac{\pi}{2})}\sspeq R\,\frac{x\,(1-x)}{1+x},\ \overline{L,D}\sspeq R\,\frac{1-x}{\tan(2\,\beta\sspm \frac{\pi}{2})}\sspeq R\,\frac{2\,x}{1+x},\ \overline{F,L}\sspeq\frac{z\,R}{\sin(2\,\beta\sspm \frac{\pi}{2})}\sspeq z\,R\,\frac{1-y}{y} \sspeq R\,\frac{(1+x^2)\, (1-x)}{2\,(1+x)},\ \overline{L,B^{\prime}}\sspeq\ R\, \frac{1-x}{\sin(2\,\beta\sspm \frac{\pi}{2})} \sspeq R\,\frac{1\sspp x^2}{1+x},\ \overline{L,C}\sspeq\ R\,\frac{v}{\tan(2\,\beta\sspm \frac{\pi}{2})} \sspeq R\,\frac{1-x}{1+x}, \  \overline{F,P}\sspeq \frac{p\,R}{\cos\,\beta }\sspm \overline{V,F} \sspeq R\,\frac{\sqrt{1+x^2}\,(2-x)}{2}, \ \overline{C,P^{\prime}}\sspeq R\,(1\sspm b\, \cos\,\beta)\sspeq R\, \frac{2-x}{2}}. \pn
With the origin $O\sspeq B$ in the $(\hat x,\hat y)$-plane the straight lines $g$ and  $g'$ are given by:\pn
\dstyle{g:\ \hat y \sspeq -x\,(\hat x \sspp R\,(1\sspp v))\sspeq -x\,\hat{x}\sspm R\frac{1+x^2}{2} ,\ \ g':\ \hat y\sspeq \frac{\hat x}{x}}.
\pbn
{\bf Observation $\bf 1^{\prime}$:} If $x$ varies from $0$ to $R$ then $P$ moves on the middle line \dstyle{\hat y \sspeq \frac{R}{2}} from $\hat x \sspeq 0$ to \dstyle{\frac{R}{2}}.
\psn
Like depicted in {\it Figure 3}, one now folds the corners $B$ {\sl and} $C$ onto the same point $C^{\prime}\sspeq B^{\prime}$ on the opposite side. The distance of $C^{\prime}$  from corner $A$ is $x\,R$.
\psn
\parbox{16cm}{
{\includegraphics[height=8cm,width=.5\linewidth]{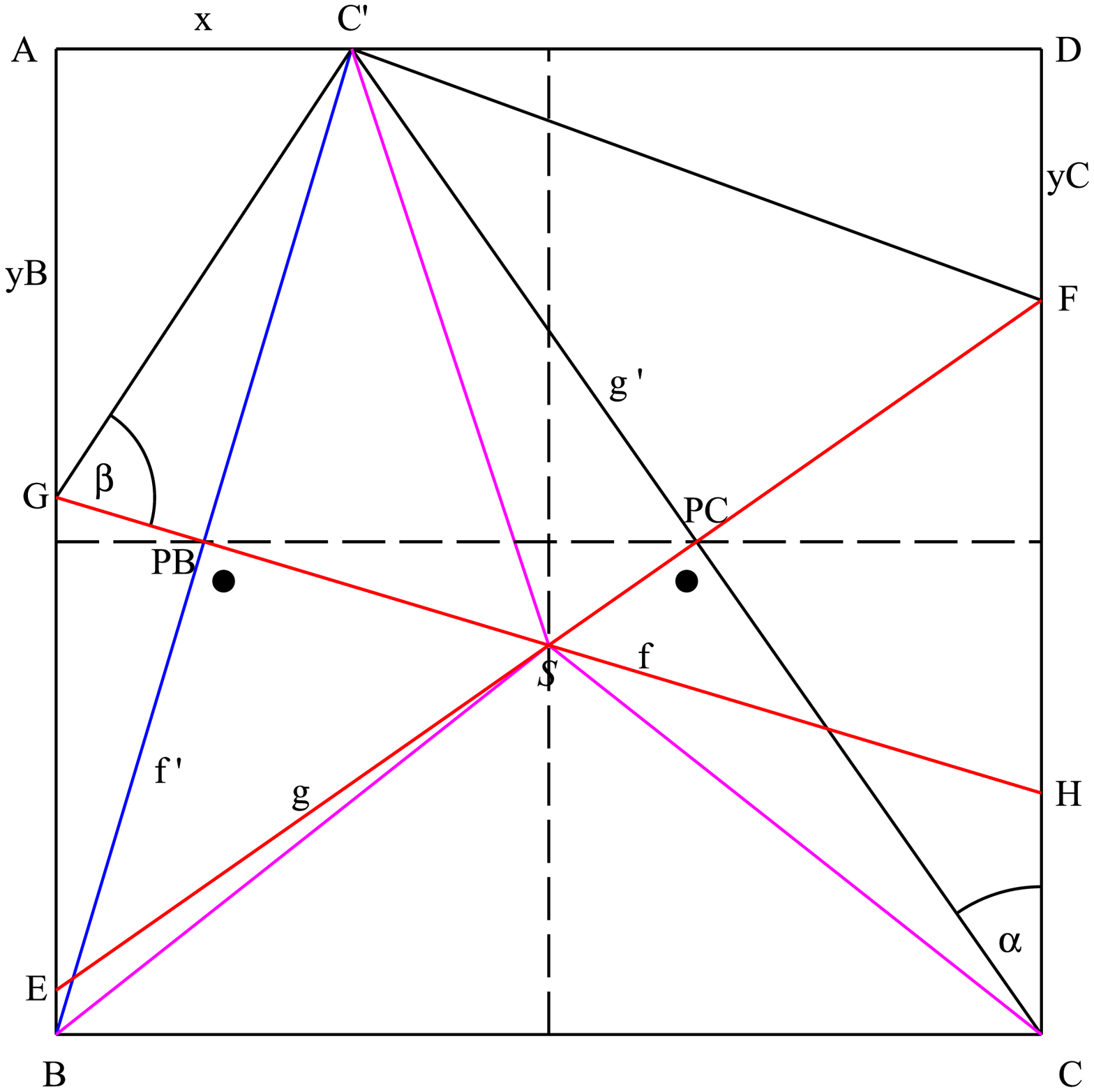}}
{\includegraphics[height=8cm,width=.5\linewidth]{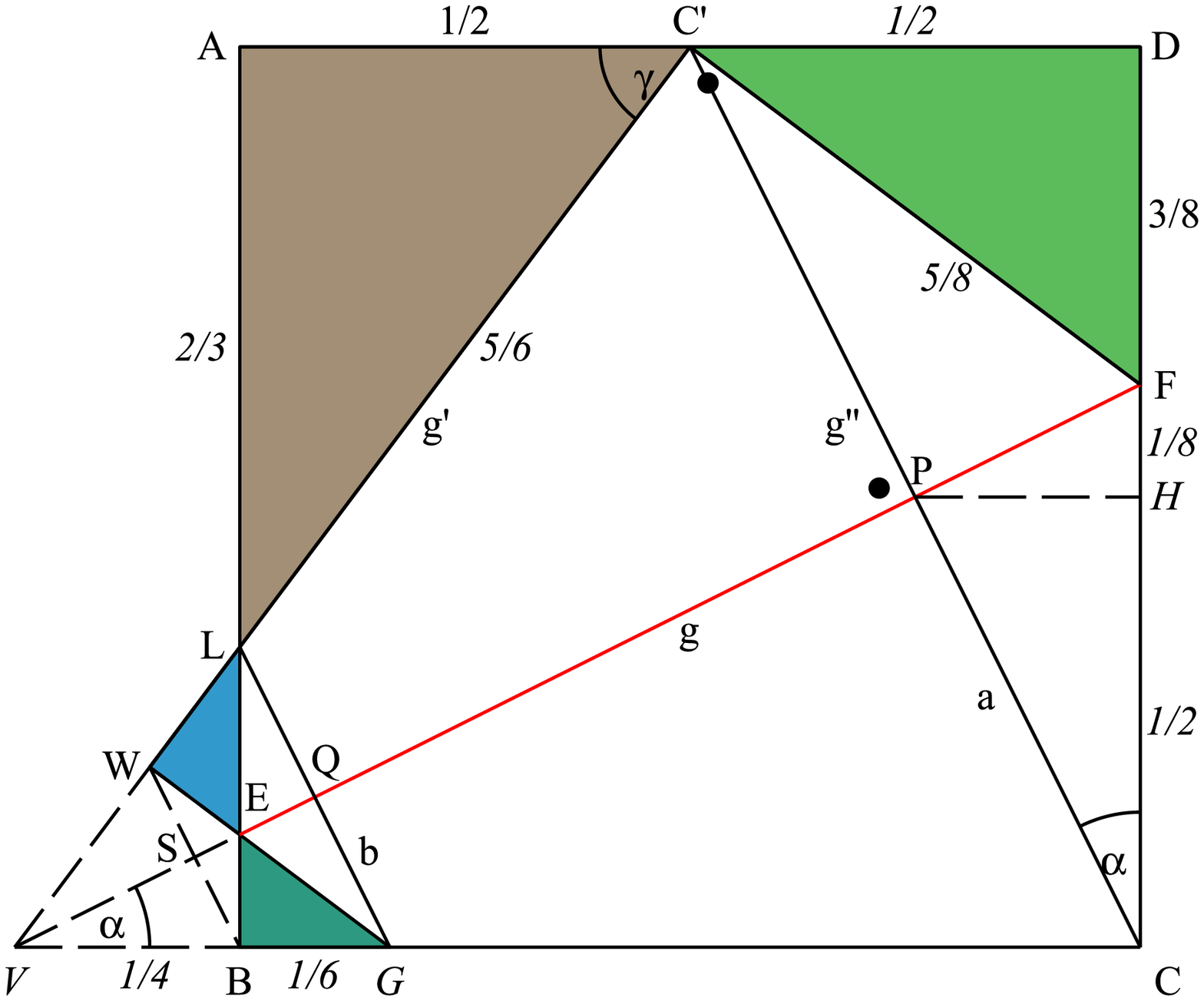}}
}
\psn
\hskip .5cm Figure 3:\  Folding $C$ and $B$ onto $C^\prime\sspeq B^{\prime}$ \hskip 1.5cm Figure 4:\  {\sl Haga}'s triple of Egyptian triangles \pbn
\vfill 
\eject
\noindent
{\bf Figure 3}:\psn
$\overline{A,B}\sspeq  R \sspeq  \overline{A,D},\ \overline{A ,C^{\prime}}\sspeq R\, x,\  \overline{D ,F}\sspeq R\, y_C,\ \overline{A ,G}\sspeq R\, y_B,\  \overline{PC ,C^{\prime}}\sspeq \overline{PC ,C}\sspfed   R\, c,\  \overline{PB ,C^{\prime}}\sspeq \overline{PB ,B}\sspfed   R\, b,\  \angle(E,PC,C) \sspeq \frac{\pi}{2} \sspeq \angle(B,PB,H),\  \angle(D,C,C^{\prime})  \sspeq \alpha \sspeq \angle(PC,C^{\prime},F),\  \angle(PB,G,C^{\prime})  \sspeq \beta \sspeq \angle(B,G,PB)$, \ $\angle(D,F,C^{\prime})\sspeq 2\,\alpha$, \ \dstyle{\angle(A,C^{\prime},G)\sspeq 2\,\beta\sspm \frac{\pi}{2}}, \ \dstyle{\angle(B,C^{\prime},C)\sspeq \alpha\sspm \beta\sspp \frac{\pi}{2}}\,.\psn
The two creases (in red) are $\overline{E,F}$, on the straight line $g$, and $\overline{G,H}$, on the straight line $f$.\psn 
$PC$ is the intersection point of the perpendicular straight lines $g$  and $g'$ connecting the points $C^{\prime}$ and $C$.  $PB$ is the intersection point of the perpendicular straight lines $f$  and $f'$ connecting the points $C^{\prime}$ and $B$. $S$ is the intersection point of the two creases $g$ and $f$.\psn
The analytic data depending on $R$ (usually taken as $1$ length unit) and $x$, is (see the above data for {\it Figures 1} and {\it 2}) :\pn 
\dstyle{ 2\,c\sspeq \sqrt{1\sspp (1-x)^2},\ 2\,b\sspeq \sqrt{1\sspp x^2},\   \overline{F,C^{\prime}}\sspeq \overline{F,C}\sspeq R\,(1\sspm y_C)\sspeq R\,\frac{2-2\,x+x^2}{2},\  \overline{G,C^{\prime}}\sspeq \overline{G,B}\sspeq R\,(1\sspm y_B) \sspeq R\,\frac{1+x^2}{2}.}.\psn
In the $(\hat x,\hat y)$-plane with origin $O\sspeq B$ the coordinates of $PC$ and $PB$ are \dstyle{\left[R\,\frac{1+x}{2} ,R\,\frac{1}{2}\right]} and  \dstyle{\left[R\,\frac{x}{2} ,R\,\frac{1}{2}\right]}, respectively.\psn
The straight lines of the two creases $g$ and $f$ are \dstyle{\hat y\sspeq (1-x)\,\left (\hat x\sspp R\,\frac{x^2}{2\,(1-x)}\right )} and \pn
\dstyle{\hat y\sspeq -x\,\left (\hat x\sspp R\,\frac{1\sspp x^2}{2\,x}\right)}, respectively. This leads to the intersection point $S$ with coordinates \dstyle{\left[\frac{R}{2}, \frac{R}{2}\,(1\sspm x\sspp x^2)\right]}.\pbn
Because $S$ lies on the vertical median of the square one obtains \dstyle{\overline{B,S} \sspeq \overline{C,S}\sspeq \frac{R}{2}\,\sqrt{1+(1-x+x^2)^2}\sspeq \frac{R}{2}\,\sqrt{2-x\,(1-x)\,(2-x\,(1-x))} \sspeq \frac{R}{2}\, \sqrt{(1\sspp (1-x)^2)\,(1\sspp x^2)}}. From crease $g$ it is clear that  $\overline{C,S}\sspeq  \overline{C^{\prime},S}$.
This proves analytically the following second theorem of  {\sl  Kazuo\ Haga}  (named like this in \cite{Bellos}, p. 115-116).\psn
{\bf Theorem 1} ({\sl K. Haga}\ see {\it Figure 3})  \psn
{\bf 1)} The intersection point $S$ of the two creases $\overline{E,F}$ and $\overline{G,H}$ lies on the vertical median of the square of length $R$, for each $x\sspin [0,R]$. \psn
{\bf 2)} The length of the three lines $\overline{C^{\prime},S}$, $\overline{B,S}$  and 
$\overline{C,S}$ are identical, namely $R\,\sqrt{2-x\,(1-x)\,(2-x\,(1-x))} $.
Moreover,\pn
{\bf 3)} The intersection points $PC$ and $PB$ lie on the horizontal median of the square of length $R$, and their distance is \dstyle{\frac{R}{2}}, independently of $x\sspin [0,R]$.\pn
The first part of {\bf 2)} follows immediately from the $x$-coordinate of $S$, {\ie} from {\bf 1)}. For the second part also the $y-$coordinate of $S$ is needed.
\pbn
{\bf Problem 2:\ Haga's Theorem on Egyptian Triangles}\psn
 This is found in the book of {\sl Bellos} \cite{Bellos} on p. 114, and has also been attributed to {\sl Kazuo Haga}. Three Egyptian triangles (or scaled Pythagorean triangles) appear when folding a corner (vertex) of a square sheet of paper (here $C$) onto the midpoint of one of the non-adjacent sides (see {\it Figure 4} point  $C^{\prime}$). The crease is the straight line $g$. The right triangles are  $T_1 \sspeq \triangle(F,D,C^{\prime})$, $T_2 \sspeq \triangle(C^{\prime},A,L)$  and  $T_3 \sspeq \triangle(L,W,E)$. $T_3^{\prime}\sspeq \triangle(G,B,E)$ folds onto $T_3$. Each of these right triangles has rational side lengths, and they are scaled $(3,4,5)-$Pythagorean triangles. 
\vfill
\eject
\noindent
{\bf  Theorem 2} ({\sl K. Haga})  \psn 
If the length of the square $R$ is taken as $1$ length unit then the sides of the three right triangles of {\sl Figure 4} have side lengths: \psn
 \dstyle{T_1:\, \left(\frac{3}{8},\, \frac{1}{2},\,\frac{5}{8}\right), \ \ T_2:\, \left(\frac{1}{2},\, \frac{2}{3},\,\frac{5}{6}\right), \ \ T_3:\, \left(\frac{1}{8},\, \frac{1}{6},\,\frac{5}{24}\right)}.\psn
Therefore, if the length of the square $R$ is chosen as $24$ length units, these triangles become {\sl Pythagorean} triangles  $\widehat{T_1}:\, (9,\,12,\,15)\sspeq 3*(3,\,4,\,5)$, $\widehat{T_2}:\, (12,\,16 ,\,20)\sspeq 4*(3,\,4,\,5)$ and  $\widehat{T_3}:\, (3,\,4 ,\,5)$.\psn
{\bf Proof}: The notation, with the length of the side of the square being $R$ length units, is: $C$ maps to $C^{\prime}$, $B$ maps to $W$. $L$ is the intersection point of $\overline{V,C^{\prime}}$ (the straight line $g^{\prime}$ with segment $\overline{W,L}$) and $\overline{A,B}$. $L$ maps to $G$.  $P$ is the midpoint of $\overline{C^{\prime},C}$,  $Q$ is the midpoint of $\overline{L,G}$ and $S$ is the midpoint of $\overline{W,B}$. The angle $\gamma\sspeq \angle(L,C^{\prime},A)$ equals \dstyle{2\,\alpha} because $\angle(E,L,W) \sspeq \frac{\pi}{2}\sspm 2\,\alpha \sspeq \angle(C^{\prime},A,L)$ \pn 
Similar right triangles with angle $\alpha$ are $\triangle (V,B,E)$, $\triangle (V,C,F)$, $\triangle (F,P,C)$, $\triangle (B,S,V)$,  $\triangle (B,E,S)$, $\triangle (L,E,Q)$, and their mirrors obtained by folding along $g$. The four shaded right triangles with angle $\gamma\sspeq 2\,\alpha$ are also similar. \dstyle{\tan\,\alpha\sspeq \frac{1}{2},\ \sin\,\alpha\sspeq \frac{\tan\,\alpha}{\sqrt{1+(\tan\,\alpha)^2}} \sspeq \frac{1}{5}\,\sqrt{5},\ \cos\,\alpha\sspeq  \frac{1}{\sqrt{1+(\tan\,\alpha)^2}} \sspeq  \frac{2}{5}\,\sqrt{5}},\  \dstyle{\tan(2\,\alpha)\sspeq \frac{2\,\tan\,\alpha}{1\sspm (\tan\,\alpha)^2} \sspeq \frac{4}{3},\  \sin(2\,\alpha)\sspeq  \frac{4}{5},\  \cos(2\,\alpha)\sspeq  \frac{3}{5}\, .}  \pn
The analytic data is:  \dstyle{ \overline{A,C^{\prime}} \sspeq \frac{1}{2}\,R\sspeq \overline{C^{\prime},D},\ \overline{C^{\prime},P}\sspeq a\sspeq  \overline{P,C} \sspeq  R\,\frac{\sqrt{5}}{4},\ \overline{L,Q}\sspeq b\sspeq \overline{Q,G}\sspeq  R\,\frac{\sqrt{5}}{12},\ \overline{W,E}\sspeq \overline{E,B} \sspeq R\,\frac{1}{8}, \ \overline{W,L}\sspeq \overline{B,G} \sspeq R\, \frac{1}{6},\ \overline{L,E}\sspeq \overline{E,G}\sspeq R\, \frac{5}{24},\ \overline{V,B}\sspeq \overline{V,W}\sspeq R\,\frac{1}{4},\  \overline{V,F}\sspeq R\,\frac{5}{8}\,\sqrt{5},\  \overline{V,E}\sspeq R\,\frac{1}{8}\,\sqrt{5},\ \overline{E,Q}\sspeq R\,\frac{1}{24}\,\sqrt{5},\ \overline{Q,P}\sspeq R\,\frac{1}{3}\,\sqrt{5},\ \overline{V,S}\sspeq R\,\frac{1}{10}\,\sqrt{5},\ \overline{S,B}\sspeq R\,\frac{1}{20}\,\sqrt{5},\ \overline{S,E}\sspeq R\,\frac{1}{40}\,\sqrt{5}}. \hskip 12.5cm $\square$
\pbn
{\bf Problem 3: Solving Third 0rder Equations using Origami}\psn
It is well known that cubic (and higher) order equations cannot be solved  geometrically using only an (unmarked) ruler and a compass (see e.g., Wantzel\cite{Wantzel}, Adler \cite{Adler}. \S 36, pp. 188-195). It is also known that the general cubic equation can be solved geometrically with two right angles (at least one of them should have a scale in order to mark the absolute value of the coefficients) \cite{Klein}, p. 267, Abb. 150, adapted from \cite{Adler}, pp. 259-261, Fig. 156 (where we use the solving chain of straight lines $A,\,X,\,Y,\,E$ with $\angle(B,A,X)\sspeq \hat \omega$ (not the $\omega$ of the figure), \dstyle{x\sspeq \tan\,\hat\omega\sspeq \frac{\overline{B,X}}{a_0} \sspeq \frac{\overline{Y,C}}{\overline{X,C}} \sspeq \frac{\overline{D,E}}{\overline{Y,D}}} and \dstyle{\overline{X,C}\sspeq |a_1|\sspm \overline{B,X}\sspeq |a_1|\sspm a_0\,x}, \dstyle{\overline{D,Y}\sspeq |a_2| \sspp  \overline{Y,C}\sspeq |a_2| \sspp x\, \overline{X,C}}. With  \dstyle{\overline{D,E}\sspeq a_3\sspeq x\, \overline{D,Y}} this leads to the cubic equation $x^3\sspm |a_1|\, x^2\sspm |a_2|\,x\sspp  a_3 \sspeq 0$ if one takes $a_0\sspeq 1$. Note that  $a_0$ and $a_1$ are supposed to have opposite signs because, coming from $\overline{A,B}$, one takes a $90^o$ left turn at $B$ to get to $C$. Similarly,  $a_1$ and  $a_2$ have like signs  because, coming from $\overline{B,C}$, one takes a left turn at $C$ to get to $D$. Then $a_2$ and $a_3$ have again opposite signs because of the left turn at $D$. The length of the lines are always positive. These sign rules are taken from \cite{Klein}, and in \cite{Adler} a different solving path, namely $A,\,F,\,G,\,H$, has been chosen. This type of figure is also found in \cite{Huzita2}, p. 198, referring to  {\it Fig. 3} on p. 207 (where the top vertex is $A^{\prime}$ which is connected to $B^{\prime}$ on line $B$. The distance between $A^{\prime}$ and $I$ is $x$, and the distance between $A$ and $B^{\prime}$ is $y$). \psn
\vfill
\eject
\psn
\parbox{16cm}{
{\includegraphics[height=8cm,width=.5\linewidth]{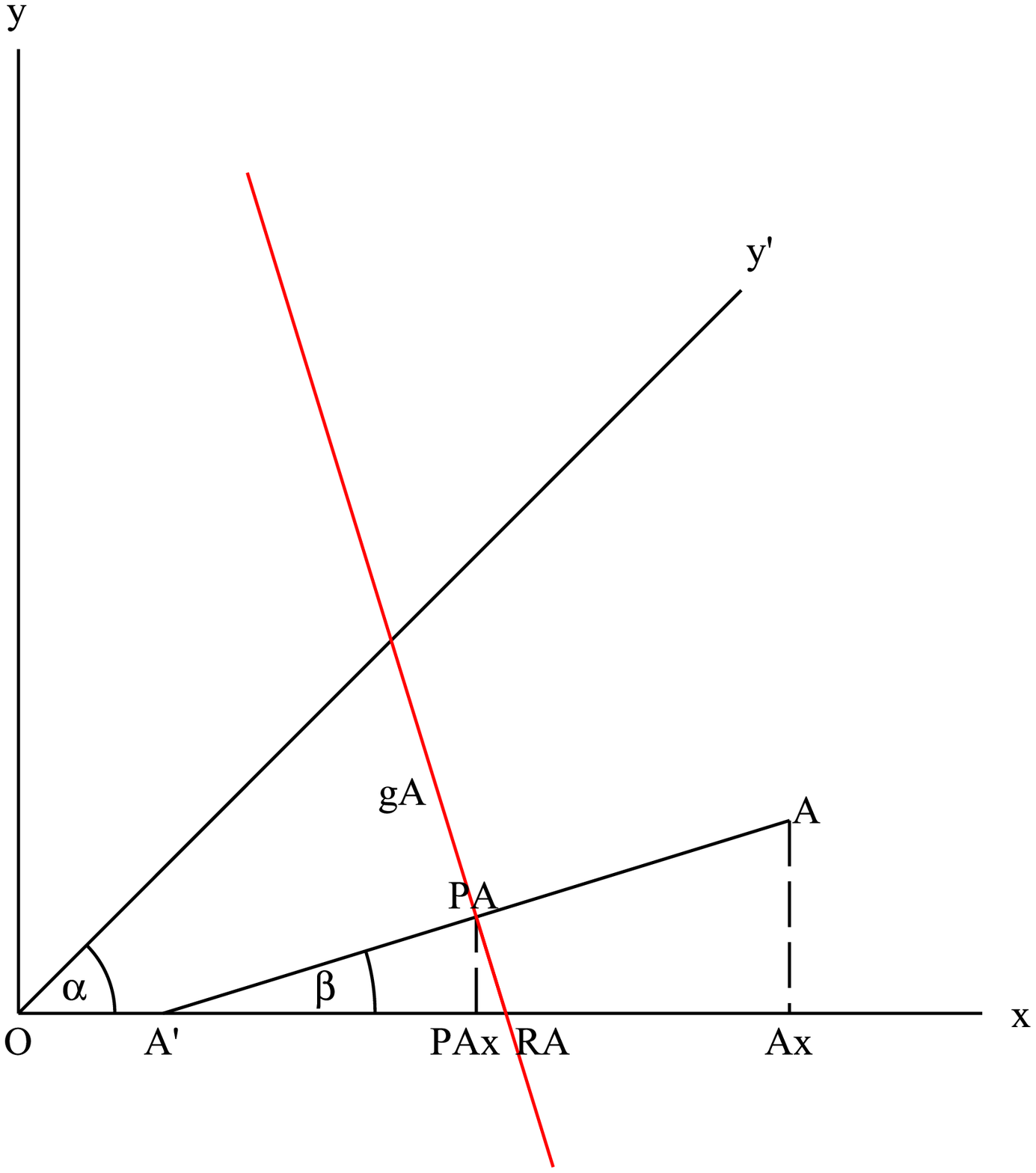}}
{\includegraphics[height=8cm,width=.5\linewidth]{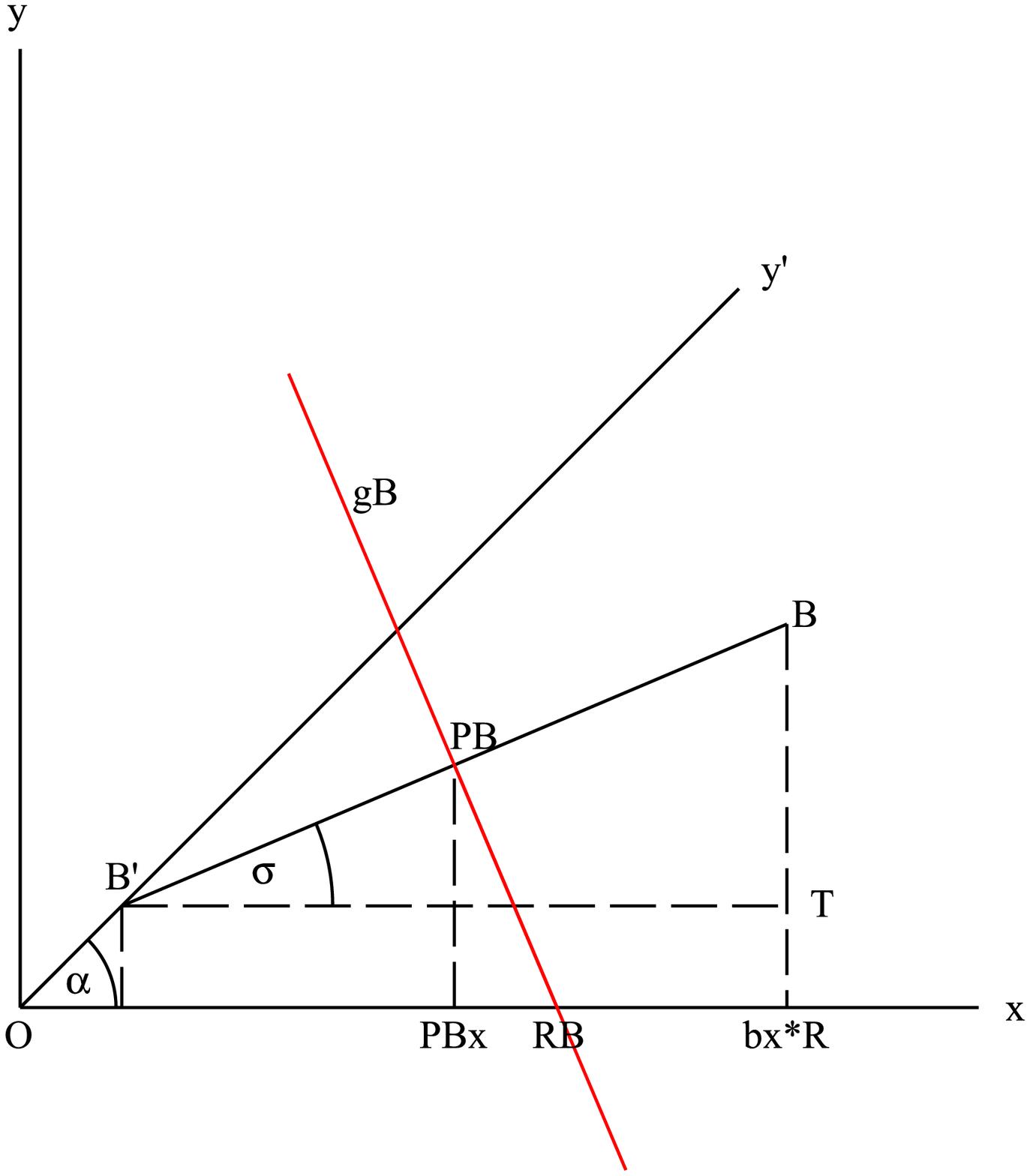}}
}
\psn
\hskip 2.5cm Figure 5:\  Folding $A$ onto $A^\prime$ \hskip 3cm Figure 6:\  Folding $B$ onto $B^\prime$  \pbn
In order to solve the general cubic equation ${\bf X^3\sspp a\,X^2\sspp b\,X\sspm c\sspeq 0}$ with origami, following Huzita \cite{Huzita2}, p. 197 and Fig. 2, one first folds , like in {\sl Figure 5}, a point $A$ onto $A^{\prime}$ on the $x$-axis. The coordinates in the $[x,y]$-plane with length unit $R$ are $A:\, [a_x\,R,\,a_y\, R]$ and $A^{\prime}:\, [\tilde x\,R,\,0]$. The $y^{\prime}$ axis with angle $\alpha$ has been added for later purposes, and is not relevant for this folding. $g_A$ is the crease perpendicular to  $\overline{A^{\prime},A}$ with intersection point $P_A:\, [P_{A,x},P_{A,y}]$. The crease hits the $x$-axis at point $R_A$. The inputs are $a_x,\, a_y$ and $\tilde x\sspeq \frac{\overline{O,A^{\prime}}}{R}$. Instead of $\tilde x$ we shall use $X\sspdef a_x\sspm \tilde x$. The right triangles $\triangle (A^{\prime},A_x,A)$ and $\triangle\, (P_A,P_{A,x},R_A)$ are similar. The following analytic expressions are found immediately.\psn
\dstyle{\overline{A,P_A}/R \sspfed a\sspeq \overline{P_A,A^{\prime}}/R\sspeq \frac{1}{2}\,\frac{a_y}{\sin\,\beta} \sspeq  \frac{1}{2}\,\frac{X}{\cos\,\beta} }.\  \dstyle{\tan\,\beta\sspeq \frac{a_y}{X} }.\ \dstyle{\overline{A^{\prime},P_{A,x}}  \sspeq \overline{P_{A,x}, A_x} },\  \dstyle{\overline{P_A,P_{A,x}}  \sspeq \frac{a_y\,R}{2}},\  \dstyle{\overline{P_{A,x},R_A }\sspeq \frac{R}{2}\, a_y\,tan\,\beta},\ \dstyle{\overline{O,R_A}\sspeq \left( a_x\sspm X \sspp \frac{a_y}{2\,tan\,\beta} + \frac{1}{2}\,a_y\,tan\,\beta\right)\,R\sspeq \frac{1}{2}\,\left(\frac{a_y^2}{X} \sspp 2\,a_x \sspm X\right) \, R}. The straight line $g_A$ (the crease) satisfies $\frac{y}{R}\sspeq -\frac{X}{a_y}\,\left(\frac{x}{R}\sspp\frac{1}{2}\,(X \sspm \frac{a_y^2}{X}\sspm 2\,a_x)\,\right)$.\psn
Then a similar folding, shown in {\it Figure 6}, is done, in order to map a point $B$ (different from $A$) onto a point $B^{\prime}$ on the $y^{\prime}$-axis, forming some angle $\alpha$ from \dstyle{\left(0,\frac{\pi}{2}\right)} with the $x-$axis. The crease $g_B$ is perpendicular to the line $\overline{B^{\prime},B}$ with  intersection at the midpoint $P_B$. It hits the $x-$axis at the point $R_B$. The $(x,y)$ coordinates of $B^{\prime}$ are $[b_x^{\prime}\,R,\, b_y^{\prime}\,R]$.\ $\angle(B,B^{\prime},T)\sspeq \sigma$. The inputs are $B: [b_x\,\,R,b_y\,R]$ and $\tilde y\sspeq \overline{0,B^{\prime}}/R$. The right triangles $\triangle\, (B^{\prime},T,B)$ and $\triangle\, (P_B,P_{B,x},R_B)$ are similar. One finds:
\psn
$b_y^{\prime}\speq \tilde y\, \sin\,\alpha$,\ \ $b_x^{\prime}\speq \tilde y\, \cos\,\alpha$.\  \dstyle{\tan\,\sigma \sspeq \frac{b_y\sspm b_y^{\prime}}{b_x\sspm b_x^{\prime}}}.\ \ \dstyle{\overline{O,R_B}\sspeq R\,\frac{1}{2}\left(b_x\sspp b_x^{\prime} \sspp \frac{(by\sspm b_y^{\prime})^2}{b_x\sspm b_x^{\prime}} \sspp 2\, b_y^{\prime}\,\tan\,\sigma \right)\sspeq R\,\frac{1}{2}\,\frac{b_x^2\sspp b_y^2\sspm \tilde y^2}{b_x\sspm b_x^{\prime}}}. \pn
The straight line $g_B$ (the crease) satisfies  \dstyle{y\sspeq -\frac{1}{\tan\,\sigma}\,(x \sspm \overline {O,R_B}) \sspeq  -\frac{b_x \sspm \tilde y\,\cos\,\alpha}{b_y\sspm \tilde y\, \sin\,\alpha}\, x   \sspp  R\,\frac{1}{2}\, \left(\frac{b_x^2\sspp b_y^2 \sspm \tilde y^2}{b_y\sspm \tilde y\, \sin\,\alpha} \right) }.\psn
In \cite{Huzita2} the coordinate axes $y^{\prime}$ and $x^{\prime}$ are used. The transformation between coordinates of a point $P:\ [p_x\,R,\,p_y\,R]$ and $[p_{x'}\,R,\,p_{y'}\,R]$ is  \dstyle{p_{x'}\sspeq p_x \sspm \frac{p_y}{\tan\,\alpha}} and  \dstyle{p_{y'}\sspeq \frac{p_y}{\sin\,\alpha}\sspeq \sqrt{1\sspp 1/(\tan\,\alpha)^2}\, p_y}.\psn
\psn
\parbox{16cm}{\begin{center}
{\includegraphics[height=8cm,width=.5\linewidth]{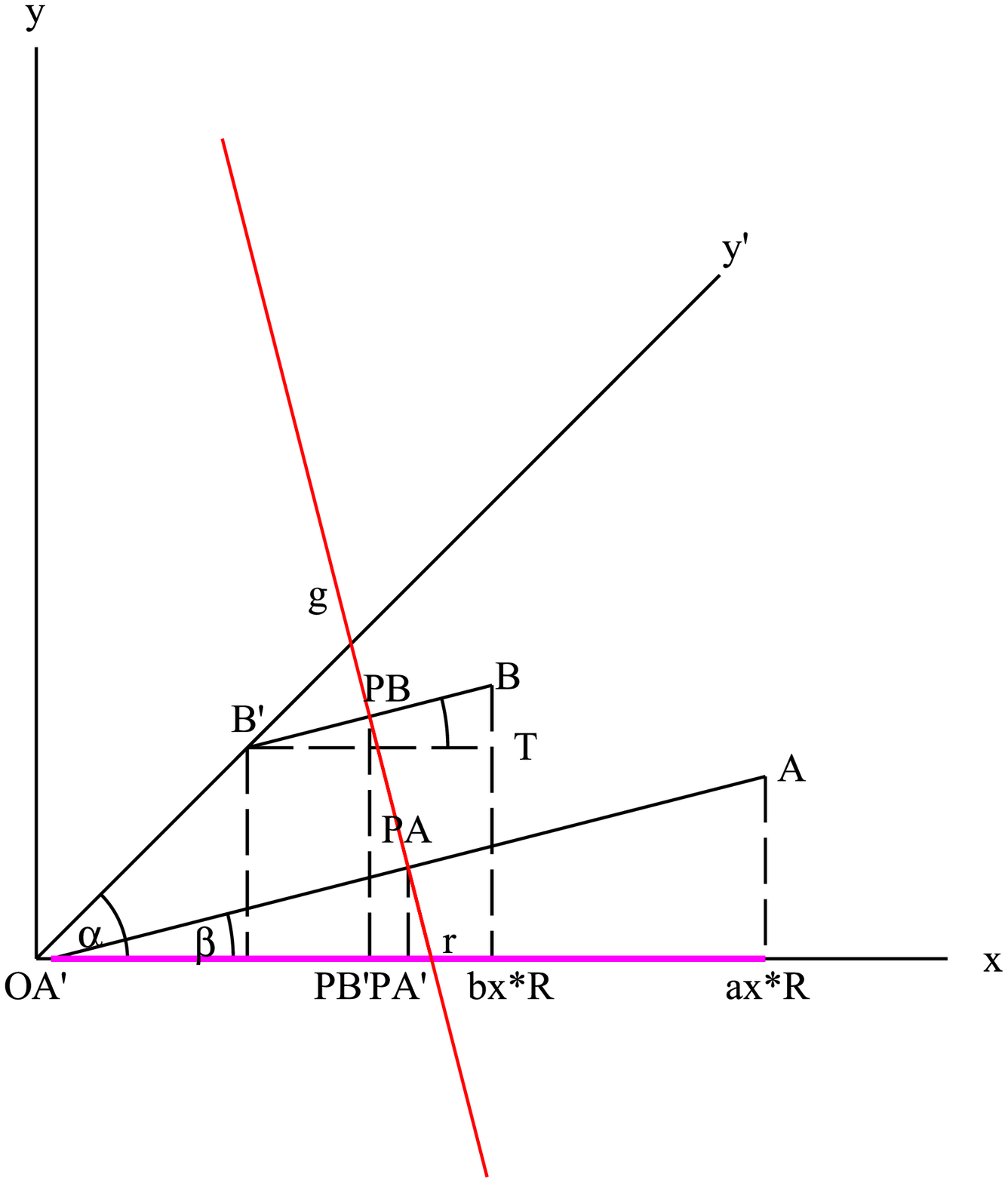}}
\end{center}
}
\psn
\hskip 2.5cm Figure 7:\  Third order eq. Folding $A$ onto $A^\prime$,  and $B$ onto $B^\prime$, $\alpha\sspeq 45^{o}$  
\pbn
The typical origami which brings at the same time one point $A$ onto $A^{\prime}$ on, say the $x$-axis, and another point $B$ onto $B^{\prime}$ on some other axis (here $y^{\prime}$)  is then shown to correspond to a cubic equation for a certain line segment (here $X\sspeq a_x\sspm a_x^{\prime}$). See {\it Figure 7} where $\alpha\sspeq 45^{o}$.\psn
{\bf Theorem 3} ({\sl Huzita} \cite{Huzita2}):  There exists a folding which brings $A:\ [a_x\,R,\,a_y\,R]$ onto $A^{\prime}:\ [a^{\prime}_x\,R,\,0] $  and $B:\ [b_x\,R,\,b_y\,R]$,onto $B^{\prime}:\ [b^{\prime}_x\,R,\, b^{\prime}_x\,\tan\,\alpha\,R]$ on the $y^{\prime}$ axis, which forms an angle \dstyle{\alpha\sspin \left(0,\frac{\pi}{2}\right)} with the $x-$axis. The solution for crease $g$ is \dstyle{y\sspeq -\frac{X}{a_y}\,(x\sspm r\,R)}, with \dstyle{r\sspeq \frac{1}{2}\,\left(\frac{a_y^2}{X}\sspp 2\,a_x -X\right)}, and $X\sspdef a_x\sspm a^{\prime}_x$ is the real solution of the cubic equation
\Beq
X^3 \sspp \left(b_x\sspm 2\,a_x\sspp \frac{b_y\sspm a_y}{\tan\,\alpha}\right)\, X^2\sspp a_y\,\left(2\,b_y\sspm a_y\sspp 2\,\frac{a_x\sspm b_x}{\tan\,\alpha}\right)\,X\sspm a_y^2\,\left(b_x\sspm \frac{a_y\sspm b_y}{\tan\,\alpha}\right)\sspeq 0\, . \nonumber
\Eeq  
Before giving the proof a remark and an example are in order.\psn
{\bf Remark:} In \cite{Huzita2} the components of $A$ and $B$ with respect to the axis $X$ and $Y$ (with $k\sspeq \cos\, \angle(X,Y)$) correspond to the above given $(a_{x'},a_{y'})$ and $(b_{x'},b_{y'})$, and $k\sspeq \cos\, \alpha$. Therefore $z$ of eq. $(1)$ on p. 197 is given by \dstyle{z \sspeq a_{x'}\sspm a^{\prime}_{x'}\sspeq  a_x\sspm \frac{a_y}{tan\,\alpha}\sspm a^{\prime}_x\sspeq X\sspm \frac{a_y}{\tan\,\alpha}} (with our $X\sspeq a_x\sspm a_{x}^{\prime}$ and $a^{\prime}_y \sspeq 0$).\psn
{\bf Example 1:} In {\sl Figure 7} we have chosen \dstyle{\alpha\sspeq \frac{\pi}{4}} and $R = 1$ length unit. With $A:\ [.8,\, .2]$ and $B:\, [.5,\, .3]$ one has the real solution of $X^3\sspm X^2\sspp 0.2\,X\sspm 0.024\sspeq 0$, which is $X\approx 0.7839279132$ (Maple 10 digits). This leads to $r\sspeq 0.4335485933$. $X$ has been indicated by the fat (magenta) line in {\it Figure 7}. $a^{\prime}_x \sspeq \tilde x \sspeq 0.8\sspm X\approx 0.0160720868$.\psn
{\bf Proof:} One combines the foldings of {\it Figure 5} and {\it Figure 6} with the constraint that the two creases, the straight lines $g_A$ and $g_B$, coincide. This leads to two equations: $(I)\ \overline{O, R_A}\sspeq \overline{O,R_B}$ and $(II)\ \beta\sspeq \sigma$.\psn
\Beqarray
&&(\overline{O,R_A}\sspm \overline{O,R_B})/R \sspeq \nonumber\\
&& \frac{1}{2\,X\,(b_x\sspm b_x^{\prime})}\,\left( X^2\,\tilde y^2 \sspp \frac{X^2\sspm 2\,a_x\,X\sspm a_y^2}{\sqrt{1\sspp (\tan\,\alpha)^2}} \,\tilde y\sspm b_x\,X^2\sspp (2\,a_x\,b_x\sspm b_x^2\sspm b_y^2)\, X\sspp a_y^2\,b_x \right)\sspeq 0\,. \nonumber
\Eeqarray
The pre-factor does not vanish and it is not divergent ($b_x\sspneq b^{\prime}_x$), therefore the bracket term has to vanish.
The other restriction is\psn
\Beq
\tan\,\sigma\sspm \tan\,\beta\sspeq \frac{b_y\sspm \tilde y\, \sin\,\alpha}{b_x\sspm \tilde y\, \cos\,\alpha} \sspm \frac{a_y}{X} \sspeq 0\ . \nonumber
\Eeq
Solving for $\tilde y$ as a function of $X$ yields
\Beq
\tilde y\sspeq \sqrt{1\sspp (\tan\,\alpha)^2}\,\frac{a_y\,b_x\sspm b_y\,X}{a_y\sspm X\,\tan\,\alpha}\ .\nonumber
\Eeq
Inserting this $\tilde y$ into the bracket term of  $(\overline{O,R_A}\sspm \overline{O,R_B})/R $  results in a factorized form given by
\Beqarray
\frac{(b_y\sspm b_x\,\tan\,\alpha)\,X}{(-a_y\sspp \tan\,\alpha\,X)^2}\,\left[(\tan\,\alpha)\,X^3+((b_x\sspm 2\, a_x)\,\tan\,\alpha \sspp (b_y\sspm a_y))\, X^2  \sspp  \right. && \nonumber \\ 
\left .\sspp ( a_y\,(2\,b_y\sspm a_y)\,\tan\,\alpha \sspp  2\,(a_x \sspm \,b_x))\,X \sspp a_y^2\,(a_y \sspm b_y \sspm b_x\, \tan\,\alpha) \right]\sspeq 0\, . \nonumber
\Eeqarray
Because $B$ does not lie on the $y^{\prime}$ axis $b_y\sspneq b_x\,\tan\,\alpha$, and $\beta \sspneq \alpha$, hence $a_y\sspneq X\,\tan\,\alpha$. Therefore, the new pre-factor does neither vanish nor diverge, and the bracket term has to vanish.\psn
Now two cases have to be considered: \psn
{\bf i)}\ $\tan\,\alpha\sspneq 0 $ and $\neq \infty$, \ie \dstyle{\alpha \sspin  \left(0,\frac{\pi}{2}\right)} and \psn
{\bf ii)}\  $\tan\,\alpha\sspeq \infty$, or $\alpha\sspeq 90^o$. \psn
The case $\alpha\sspeq 0^o$ will later be treated separately.\psn
{\bf i)}: This case leeds to the cubic equation for $X$ given in {\it Theorem 3} after dividing $\tan\,\alpha$ out. The equation for the crease $g$ is just given by $g_A$ from above (see {\it Figure 5}), with $r\sspeq \overline{O,R_A}$ and $\beta\sspeq \sigma$ from condition $(II)$. \hskip 16cm $\square$
\pbn
{\bf Special case  ii)\  $\bf{\boldsymbol \alpha}\sspeq \frac{\pi}{2}$} (see {\it Figures 8} and {\it 9})\psn
{\bf Theorem 4}  \psn
With the notation of {\it Theorem 3} and \dstyle{\alpha\sspeq \frac{\pi}{2}} the equation for $X\sspeq a_x\sspm a_x^{\prime}$ is 
\Beq 
X^3\sspp (b_x\sspm 2\,a_x)\,X^2\sspp a_y\,(2\,b_y\sspm a_y)\,X\sspm a_y^2\,b_x\sspeq 0\  \nonumber.
\Eeq  
{\bf Proof:} Extract $\tan\,\alpha$ from the bracket term of$(\overline{O,R_A}\sspm \overline{O,R_B})/R\sspeq 0$, factorized above, and observe that the original  pre-factor \dstyle{\frac{1}{2\,X\,(b_x\sspm b_x^{\prime})}} when multiplied with the factor \dstyle{\frac{(b_y\sspm b_x\,\tan\,\alpha)\,X}{(-a_y\sspp \tan\,\alpha\,X)^2}} and after extraction of the $\tan\,\alpha$ from the bracket term becomes, in the limit $\tan\alpha \sspto \infty$, \dstyle{\frac{1}{2\,X^2}}, provided $b_x\sspneq 0$. This new factor does neither vanish nor diverge, and from the bracket term the claimed cubic equation for $X$ is obtained. \pn
This result can also be reached  in the limit $\tan\alpha\sspto \infty$ from the above equation for $X$ in {\it Theorem 3}, which, however,  has been derived assuming $\tan\,\alpha\sspneq 0.$  \hskip 8cm $\square$\psn
For the case \dstyle{\alpha\sspeq \frac{\pi}{2}} see {\it Figures 8} and {\it 9}.
The other quantities are found by first folding $B$ onto $B^{\prime}$ on the $y-$axis. See {\it Figure 8}. The free parameter is $\tilde y\sspeq \overline{O,B^{\prime}}$. \dstyle{\delta\sspeq \frac{\pi}{2}\sspm \gamma,\ b\sspeq \overline{B,PB}\sspeq \overline{PB,B^{\prime}}\sspeq \frac{1}{2}\,\sqrt{(\tilde y\sspm b_y)^2\sspp b_x^2},\  \tan\,\gamma\sspeq \frac{\overline{O,PB_x}\sspm \overline{O,R_B}}{b_y\sspp b\,\sin\, \gamma},\ PB_x:\, [R\,(b_x\sspm \cos\, \gamma\, b),\,0],\   PB_y:\, [0,\,R\,(b_y \sspp \sin\,\gamma\, b)],\ RB:\, [PB_x[1]\sspm \tan\, \gamma\,PB_y[2],\, 0]}. The crease is $g:\ y\sspeq \frac{1}{\tan\, \gamma}\,(x\sspm RB[1])$.\psn
Then also $A$ is folded onto $A^{\prime}$ on the $x-$axis. This has been treated in connection with {\it Figure 5} (were the $y^{\prime}-$axis was not important). The two constraints are $RA\sspeq RB$ and $\tan\, \delta\sspeq tan\,(\frac{\pi}{2}\sspp \beta)$, \ie $\tan\, \gamma \sspeq -\tan\, \beta$. One obtains a real solution of the cubic equation for $X\sspeq a_x\sspm a^{\prime}_x $ indicated by the thick (magenta) line segment in {\it Figure 9}, and  \dstyle{\tilde y\sspeq -b_x \frac{a_y}{X}\sspp b_y}. \hskip 7cm $\square$  \psn
\parbox{16cm}{
{\includegraphics[height=8cm,width=.5\linewidth]{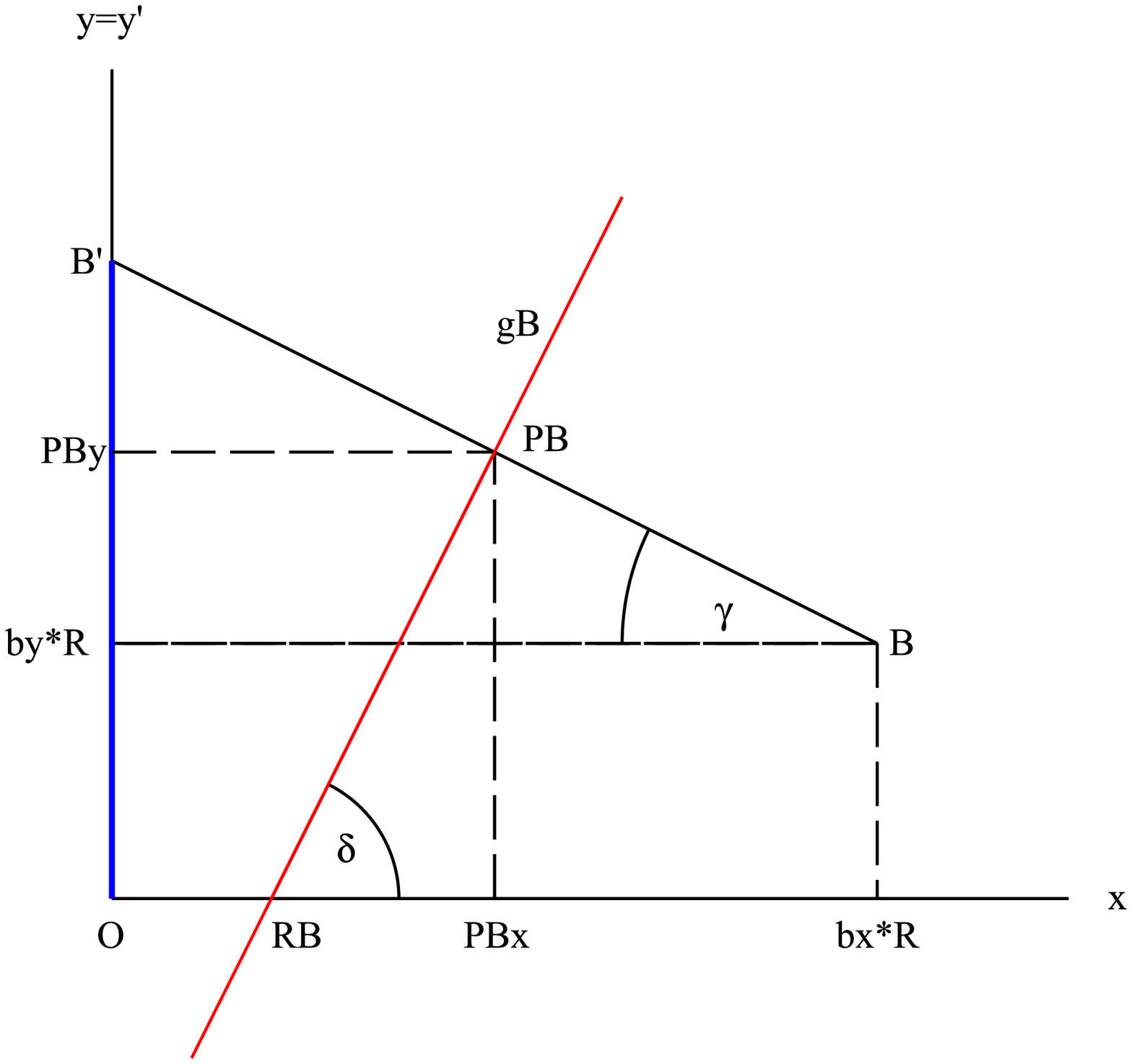}}
{\includegraphics[height=8cm,width=.5\linewidth]{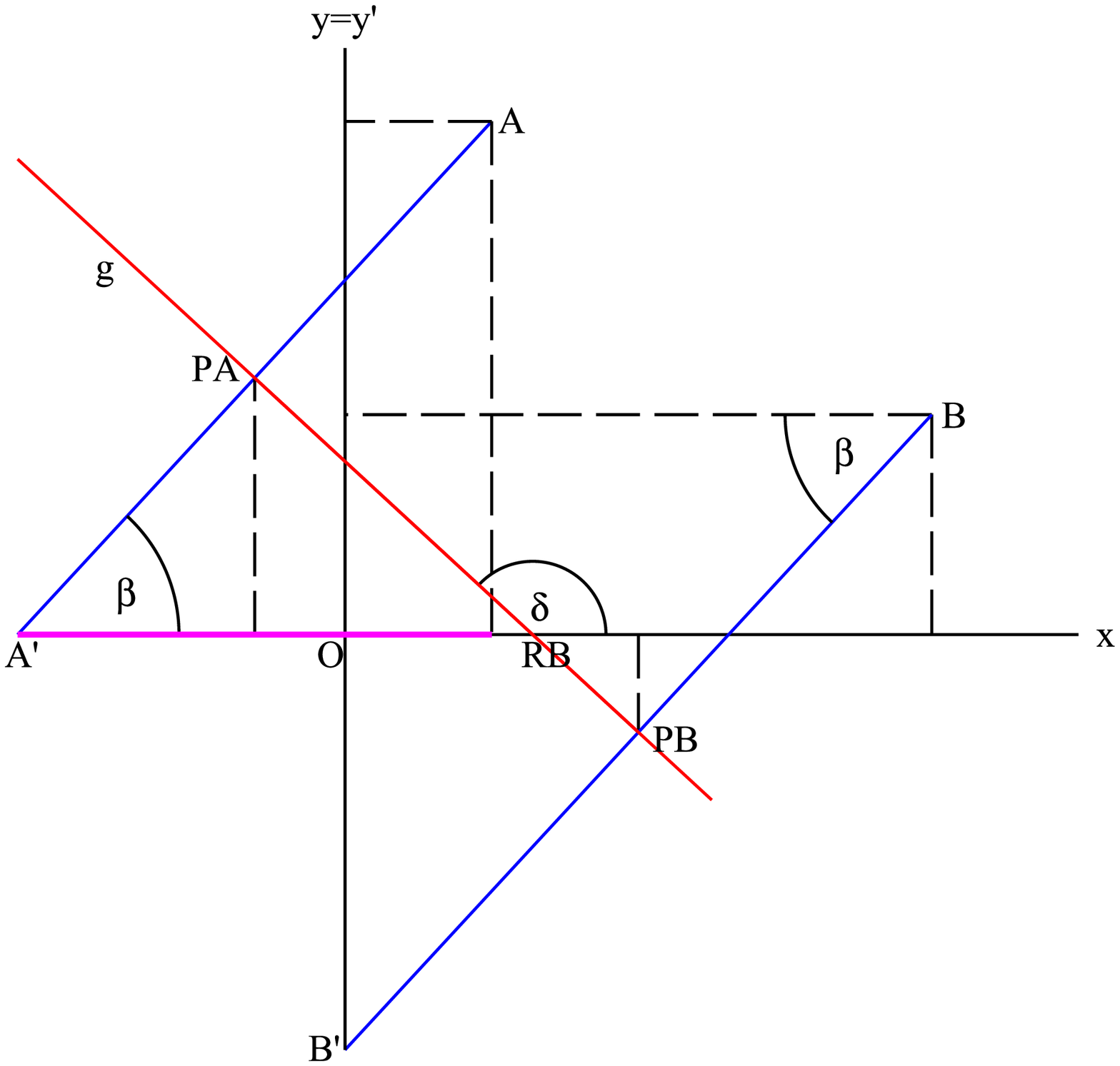}}
}
\psn
\hskip 6cm  Case $\alpha\sspeq \frac{\pi}{2}$
\psn
\hskip .5cm Figure 8:\ Folding $B$ onto $B^\prime$ \hskip 2.5cm Figure 9:\  Folding $A$ onto $A^{\prime}$ and $B$ onto $B^\prime$  \pbn
\pbn
{\bf Example 2:} $R\sspeq 1$ length unit, $A:\ [.2,.7]$, $B:\ [.8,.3]$.  $X\sspeq a_x\sspm a_x^{\prime}\approx .646416 $, $\beta\approx 47,279^o$, $\delta\sspeq (90\sspp \beta)^o\approx 137.279$, $\overline{O,R_B}\approx .2558$, $P_A:\,  [\approx -.123,.35]$, $P_B:\ [.4,\approx -.133]$ .
\pbn 
{\bf Special case $\bf{\boldsymbol \alpha}\sspeq 0$}\psn
The special case $\alpha\sspeq 0$ is obtained from folding $A$ onto $A^{\prime}$, discussed above (see {\it Figure 5}), and folding $B$ onto $B^{\prime}$ (also on the $x-$axis). One uses the formulae given above in connection with {\it Figure 5}  (with $\beta$ changed into  $\beta_A$) and replaces there $A$ by $B$ and $A^{\prime}$ by $B^{\prime}$ (with $\beta\sspeq \beta_B$). Then {\it Figure 10} is obtained by setting  $R_A\sspeq R_B$ and $\beta_A\sspeq \beta_B$.
\psn
\parbox{16cm}{\begin{center}
{\includegraphics[height=8cm,width=.5\linewidth]{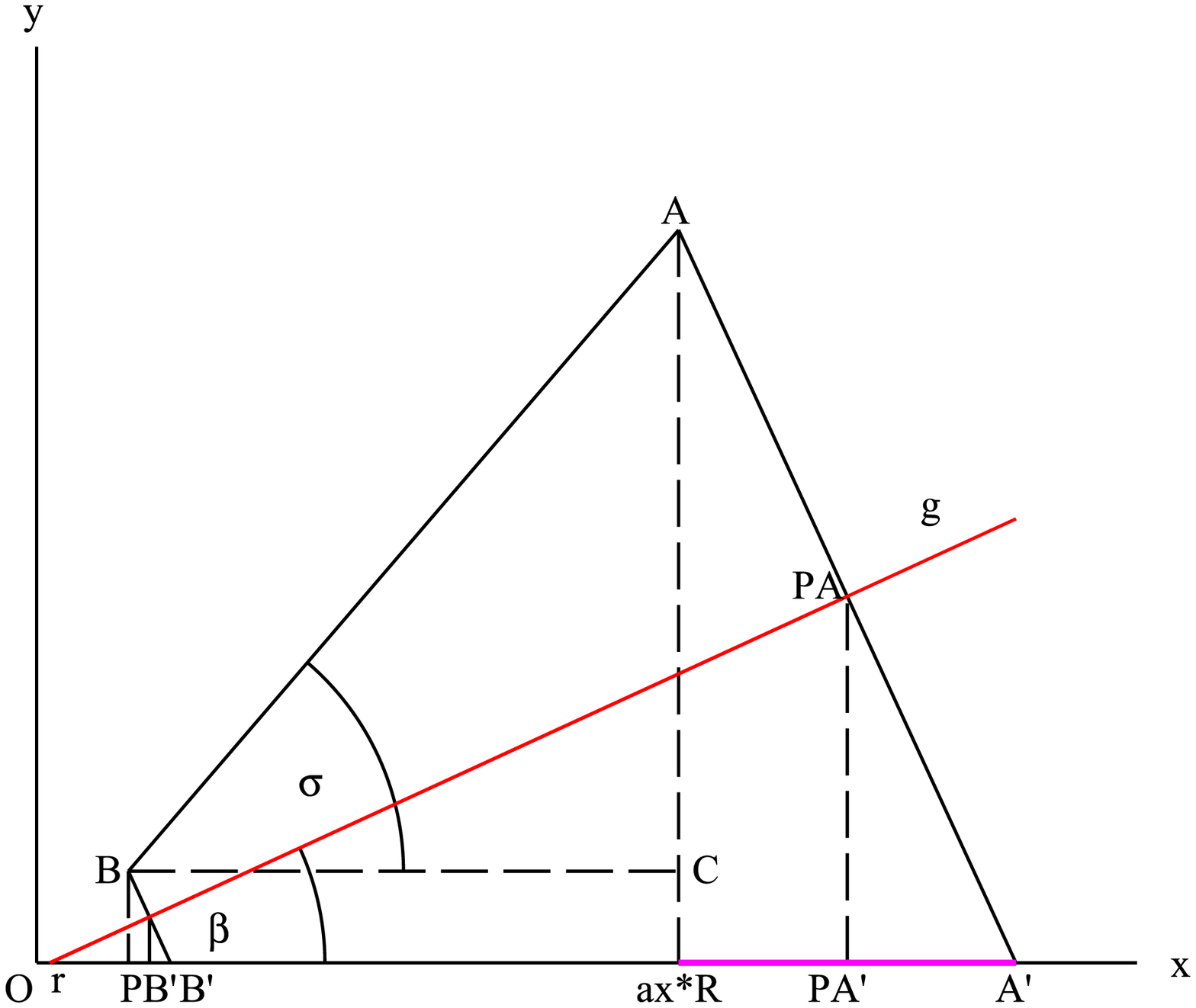}}
\end{center}
}
\psn
\hskip 3.5cm Figure 10:\ Folding $A$ on $A^\prime$  and $B$ on $B^\prime$ when $\alpha\sspeq 0$  
\vfill
\eject
\noindent
The result of identifying both creases $g_A$ and $g_B$, calling them $g$, leads to the following data:\pn
$X\sspdef a_x\sspm a_x^{\prime}$,\  $Y\sspdef b_x\sspm b^{\prime}_x$.\ Setting $\overline{O,R_A}\sspm \overline{O,R_B}\sspeq 0$ leads to the equation involving $X$ and $Y$.
\Beq
-X\sspp Y\sspp a_y^2\,\frac{1}{X} \sspm b^2_y\,\frac{1}{Y} \sspp 2\,(a_x\sspm b_x)\sspeq 0\ .\nonumber
\Eeq
The slope of $g$ is given by $\tan\,\beta_A\sspeq \tan\,\beta_B$, which expresses $Y$ in terms of X: 
\Beq
Y\sspeq \frac{b_y}{a_y}\, X \ . \nonumber
\Eeq
The mid points of  $\overline{A\,A^{\prime}} $ and $\overline{B\,B^{\prime}} $  are \dstyle{P_A:\, \left[R\,\frac{2\,a_x\sspm X}{2},R\,\frac{a_y}{2}\right]} and \dstyle{P_B:\, \left[R\,\frac{2\,b_x\sspm Y}{2},R\,\frac{b_y}{2}\right]}, respectively. The intercept is \dstyle{r\sspeq\overline{O,R_A}\sspeq \overline{O,R_B}\sspeq \frac{R}{2}\, \left(\frac{a_y^2}{X}\sspp 2\,a_x\sspm X \right)}.\pn
For the following one assumes that $0\sspneq a_y\sspneq b_y$.
Plugging $Y$ into the previous equation results in a quadratic equation for $X$:
\Beq
X^2\sspm 2\,a_y\frac{1}{\tan\,\sigma}\,X\sspm a^2_y\sspeq 0\ , \nonumber
\Eeq
with \dstyle{\tan\,\sigma\sspeq \frac{a_y\sspm b_y}{a_x\sspm b_x}}.\pn
This equation can be obtained directly from the above given analysis for non-vanishing $\alpha$ before $\tan\, \alpha$ has been divided out. Just let there $\tan\,\alpha\sspto 0$, which eliminates the $X^3$ term, and divide by $b_y\sspm a_y$. The relevant solution for $X$ is then 
\Beq
X\sspeq a_x \sspm a_{x^{\prime}}\sspeq -a_y\,\frac{1\sspm \cos\, \sigma}{\sin\, \sigma }\ .\nonumber
\Eeq
{\bf Example 3:} $R\sspeq 1$ length unit,\ $A:\, [.7,.8],\, B:\,[.1,.1]$. $\sigma\sspapprox  49.399^o$, $-X\sspapprox  0.368$, and \ $r\sspapprox 0.014$.\psn
If $a_y\sspeq b_y$, the equation for $X$  becomes linear, in fact $X\sspeq 0$.\psn 
{\bf Degenerate case: Parallel  lines with $\bf A^{\prime}$ and $\bf B^{\prime}$} \psn
As mentioned in \cite{Huzita2}, p. 197, the case of parallel lines, is also of interest. See {\it Figure 12}.
\psn
\parbox{16cm}{
{\includegraphics[height=8cm,width=.5\linewidth]{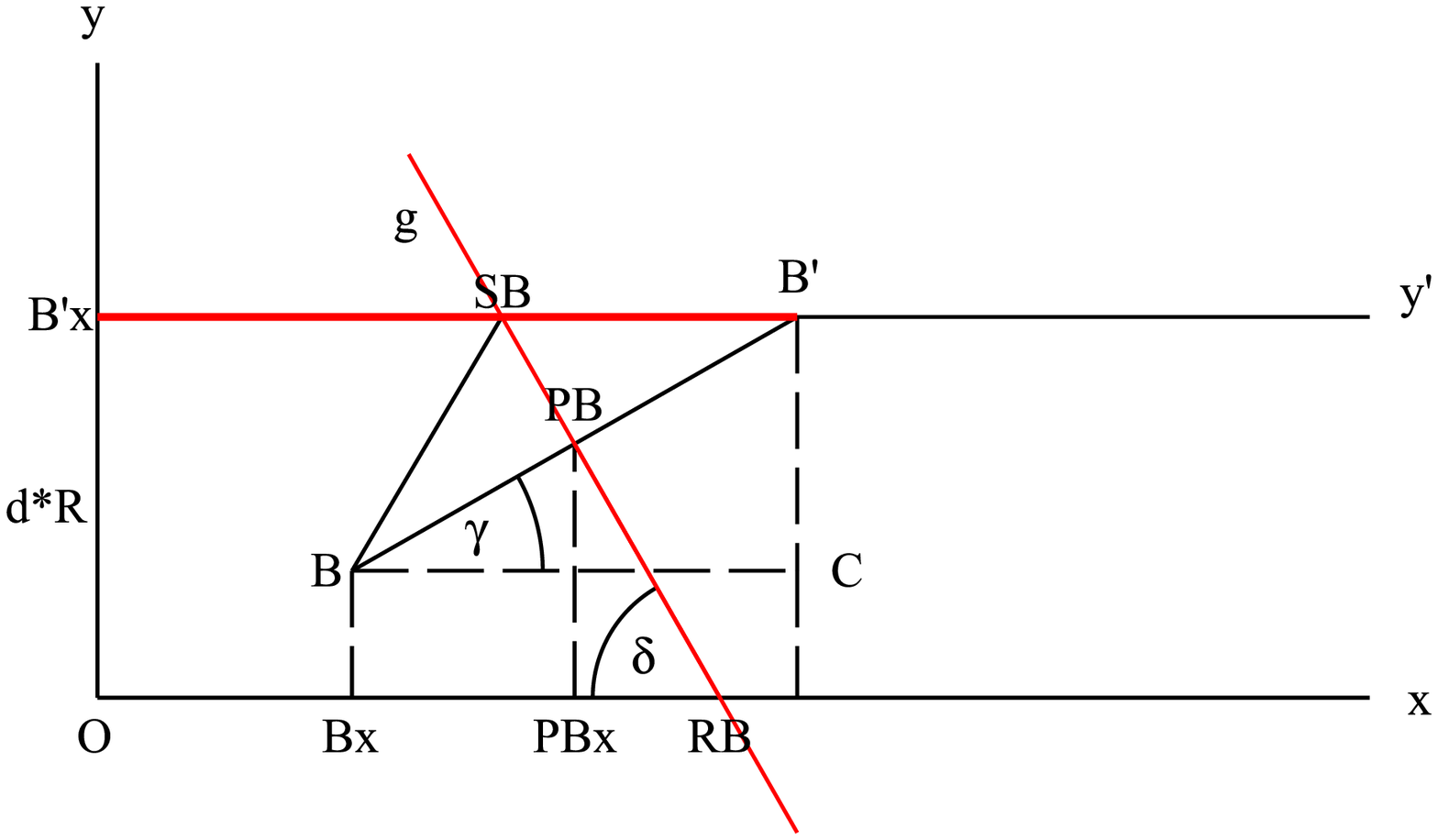}}
{\includegraphics[height=8cm,width=.5\linewidth]{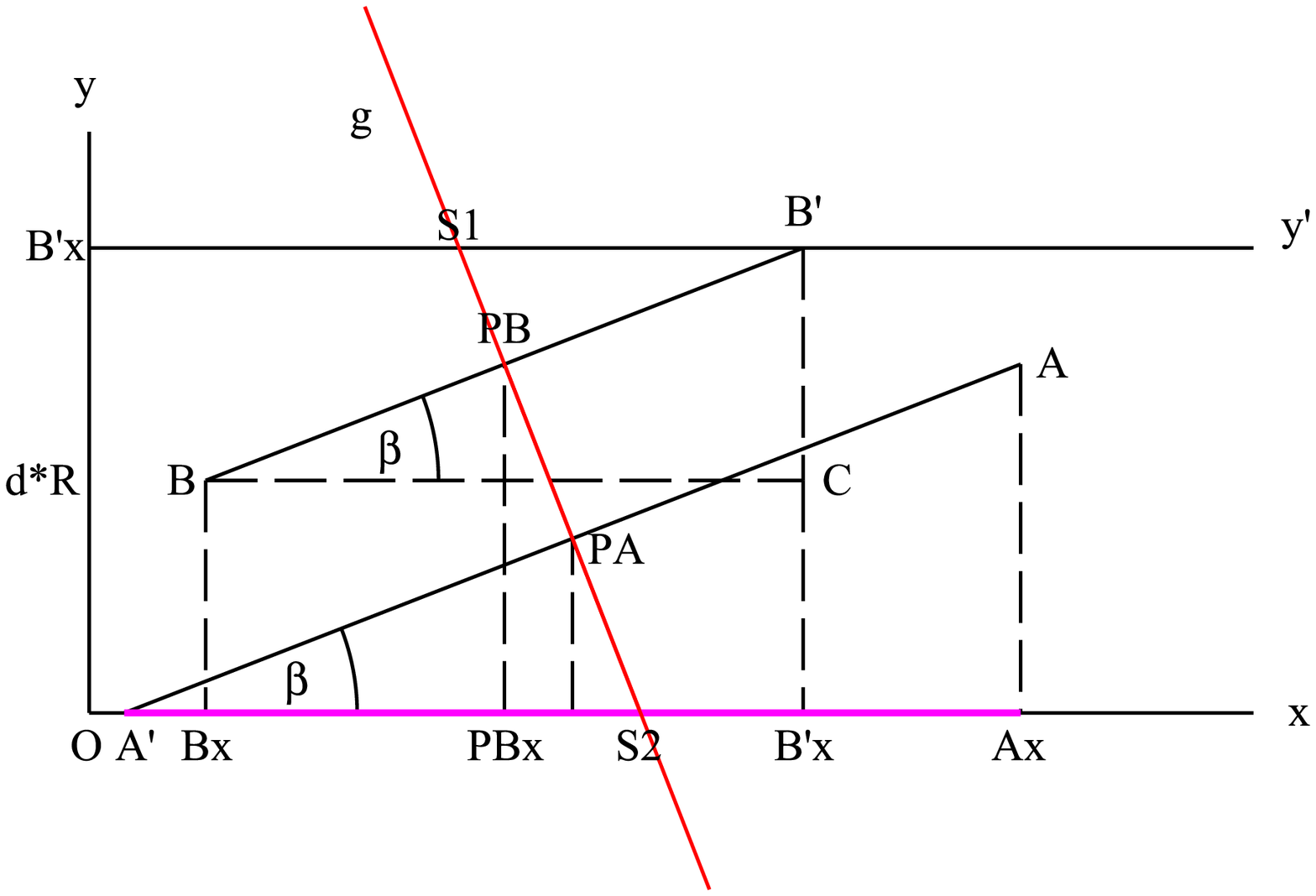}}
}
\psn
\hskip 1cm Figure 11:\  Folding $B$ onto $B^\prime$ \hskip 3cm Figure 12:\  Folding $A$ onto $A^{\prime}$ and $B$ onto $B^\prime$  \pbn
First one folds $B$ onto $B^{\prime}$ which lies on the horizontal $y^{\prime}$-axis with $y\sspeq d\,R$, hence $B^{\prime}_y\sspeq d\,R$, with some length scale $R$. This is shown in {\it Figure 11}.\psn
{\bf Figure 11:} \psn
Given $d\sspeq b^{\prime}_y$, $A$ and $B$, the free parameter is \dstyle{b^{\prime}_x\sspeq \frac{B^{\prime}_x}{R}}, the position of $B^{\prime}$ on the  $y^{\prime}$-axis. The straight line $g$ is perpendicular to the line segment $\overline{B,B^{\prime}}$ and passes through its midpoint $P_B$. This line $g$ intersects the parallel $x$- and $y^{\prime}$-axes at $R_B$ and $S_B$, respectively. The angle $\delta$ equals \dstyle{\frac{\pi}{2}\sspm \gamma}. Half the distance between $B$ and $B^{\prime}$ is $R\,b$. \pn
The formulae are: \dstyle{b\sspeq \frac{b^{\prime}_x\sspm b_x}{2\, \cos\, \gamma},\ \tan\, \gamma\sspeq \frac{b_y^{\prime}\sspm b_y}{b_x^{\prime}\sspm b_x},\ \cos\, \gamma\sspeq \frac{1}{\sqrt{1\sspp (\tan\, \gamma)^2}},\  \sin\, \gamma\sspeq \frac{\tan\,\gamma}{\sqrt{1\sspp (\tan\, \gamma)^2}} },\pn
\dstyle{PB:\ [R\,(b_x\sspp b\,\cos\,\gamma),\,R\,(b_y\sspp b\,\sin\,\gamma)],\  R_B:\ \left[R\,\frac{b^{\prime\ 2}_x\sspm b^2_x\sspp b^{\prime\ 2}_y\sspm b^2_y}{2\,(b^{\prime}_x\sspm b_x )},\, 0\right]},\pn \dstyle{g:\ y\sspeq -\frac{1}{\tan\, \gamma}\ (x\sspm  R_{B,x}),\  S_B:\ \left[R_{B,x}\sspm d\, \tan(\gamma)\,R,\, d\,R\right].} If $B$ lies on the $y^{\prime}$-axis then $\gamma\sspeq 0$ and \dstyle{\delta\sspeq \frac{\pi}{2}}, a simple special case.\psn
Now $A$ is folded onto $A^{\prime}$ on the $x$-axis at the same time as $B$ is folded onto $B^{\prime}$. See {\it Figure 12}. From the above results in connection with {\it Figure 5} one takes $R_A$ and $\tan \beta$ in terms of $X\sspdef a_x\sspm \tilde x$, where $\overline{O,A^{\prime}}\sspeq R\,\tilde x$. The solution is obtained from $R_A\sspeq R_B$ and $\tan\, \beta\sspeq \tan\,\gamma$. The latter equation can be used to eliminate $b^{\prime}_x$ by (remember that $d\sspeq b^{\prime}_y$)
\Beq
 b^{\prime}_x \sspeq \frac{1}{a_y}\, ((d\sspm b_y)\,X \sspp b_x\,a_y)\ ,\nonumber
\Eeq
provided $a_y\sspneq 0$, which we assume. If $a_y\sspeq 0$ then $b_y\sspeq b^{\prime}_y\sspeq d$, $\beta\sspeq \gamma \sspeq 0$ and \dstyle{\delta \sspeq \frac{\pi}{2}}, a simple special case. From $2\, \overline{O,R_A}\sspeq 2\,\overline{O, R_B}$  one obtains
\Beq 
\-X\sspp \frac{a_y^2}{X}\sspp 2\,a_x\sspm \frac{b^{\prime\ 2}_x\sspm b_x^2\sspp d^2\sspm b_y^2}{b^{\prime}_x\sspm b_x} \sspeq 0\ . \nonumber
\Eeq
If $b^{\prime}_x$ is inserted one finds a quadratic equation for $X$, assuming that $a_y\,X$ does neither vanish nor diverge.  Because $a_y\sspneq 0$ has been assumed, if $X\sspeq 0$ then  $a_x\sspeq a^{\prime}_x$, $b_x^{\prime}\sspeq b_x$ \dstyle{\gamma\sspeq \beta\sspeq \frac{\pi}{2}}, $\delta \sspeq 0$, another simple special case. With $a_y$ also $X$ will be finite. 
\Beq
(a_y\sspm b_y\sspp d)\,X^2\sspm 2\,a_y\,(a_x \sspm b_x)\,X\sspm a_y^2\,(a_y\sspm b_y\sspm d)\sspeq 0 .\nonumber
\Eeq
Therefore, assuming $(b_y\sspneq d \sspp a_y) $ one finds the positive solution for $X$
\Beq
X\sspeq \frac{a_y}{d\sspp a_y\sspm b_y}\,\left( a_x\sspm b_x \sspp \sqrt{(a_x\sspm b_x)^2\sspp(a_y \sspm b_y)^2\sspm d^2} \right)\ .\nonumber
\Eeq
This shows, that a solution is only possible if $d\,R$ does not exceed the distance between $A$ and $B$, which has been observed in \cite{Huzita2}. It is clear that $P_A$ and $P_B$ have coordinates which are the arithmetic mean between the corresponding coordinates of $A$ and $A^{\prime}$ and $B$ and $B^{\prime}$, e.g., \dstyle{P_{B,x}\sspeq R\,\frac{b_x\sspp b^{\prime}_x}{2}}. Note that \dstyle{\overline{B,B^{\prime}}\sspeq \frac{d\sspm b_y}{a_y}\, \overline{A,A^{\prime}}}. Therefore the trapezoid  $A^{\prime},A,B^{\prime}$ and $B$ becomes a rectangle  precisely if  $d\sspeq a_y\sspp b_y$. 
\psn
{\bf Example 4}: Put  $R\sspeq 1$ length unit and take $ d\sspeq .4,\, a_x\sspeq .8,\, a_y\sspeq .3,\, b_x\sspeq .1,\, b_y \sspeq .2 $, then  $X\sspeq a_x\sspm a^{\prime}_x\sspapprox 0.77$, $b^{\prime}_x\sspapprox 0.61$.
\pbn
\vfill 
\eject
\noindent
{\bf Application 1: The case of the heptagon equation}\psn
The minimal polynomial of the algebraic number \dstyle{\rho(7)\sspdef 2\,\cos\left(\frac{\pi}{7}\right)\approx  1.801937736} (the length ratio of the larger diagonal and the side of a regular $7$-gon)  is $C(7,x)\sspeq x^3\sspm x^2\sspm 2\,x\sspp 1$ (see \eg \cite{WLang}), {\it Table 2} and {\it section 3}). The three real zeros are known to be \dstyle{x(7;k)\sspeq 2\,\cos\left(k\,\frac{\pi}{7}\right)}, for $ k\sspeq 1,\,3,$ and $5$. They are  $x(7;1)\sspeq \rho(7)$,  $x(7;3)\approx .4450418670$ and \dstyle{x(7;5)\sspeq - 2\,\cos\left(2\,\frac{\pi}{7}\right)\approx  -1.246979604}.\psn
Here we show how these zeros are obtained by three different origamis. We also treat the standard geometric solution of this cubic equation using two right angular rulers, as explained in \cite{Klein} based on \cite{Adler} (see also {\sl von Sanden} \cite{vonSanden}, ch. III, sect. 2, pp. 55-61, with {Fig. 17} on p. 55). The corresponding {\it Figures} are {\it 13},\ {\it 14}, and {\it 15}. The slope of the $y^{\prime}$ axis is chosen as $\alpha \sspeq 90^o$, thus $y^{\prime}\sspeq y$. The monic cubic equation $C(7,x)\sspeq 0$ has sign pattern  $+,-,-,+$.  This leads, in the standard geometrical construction, to the right angle pattern $l,r,l$, with $l$ and $r$ for a $90^o$ left and right turn, respectively. One starts with some (oriented) horizontal line segment $\overline{B,C}$ of length $a_0\sspeq 1$ (for the monic case in some length unit $R$). A $90^o$ left turn gives $\overline{C,D}$ of length $a_1\sspeq 1$, then a $90^o$ right turn leads to $\overline{D,E}$ of length $a_2\sspeq 2$, and finally the $90^o$ left turn leads to $\overline{E,A}$ of length $a_3\sspeq 1$. (The starting point has been chosen as $B$ in order to comply with the later origami solution). This pattern (`Streckenzug' or line segment zig-zag) is dictated by the cubic equation and will be the same for each of the three solutions. \pn
In the origami version one needs the two perpendicular axes $y^{\prime}\sspeq y$ and $x$. As explained in \cite{Huzita2}, p. 198-199  and {\it Fig. 4} on p. 208, the $y$-axis is chosen  parallel to $\overline{CD}$, at a perpendicular distance $2\,a_0\sspeq 2$ from point $B$. The $x$-axis is parallel to  $\overline{DE}$ at a perpendicular distance $2\,a_3\sspeq 2$ from point $A$.  See the present {\it Figure 13}. In our case point $B$ has  coordinates $[-2,0]$ and $A:\, [1,2]$ (if $R\sspeq 1$) . In the standard geometrical construction of a solution to the cubic equation one has to find a point $F$ on the axis with line segment $\overline{C,D}$, here $F:\, [-1,x]$, such that a line perpendicular to $\overline{B,F}$ through $F$ hits point $G$ on the straight line with segment $\overline{D,E}$, and a perpendicular line to $\overline{F,G}$ through $G$ hits point $A$. In {\it Figure 13} the solution has $F\sspeq P_B\sspequiv PB$ and $G\sspeq P_A\sspequiv PA$. For a general cubic equation there will always be at least one real solution, and depending on its discriminant one will find either one, two or three real solutions. In general the discriminant is $Disc\sspeq p^3+q^2$, with \dstyle{q\sspdef \frac{1}{2}\,\left (2\,\frac{a_1^3}{27}\sspm\frac{a_1\,a_2}{3}\sspp a_3\right)} and \dstyle{p\sspdef \frac{1}{9}\,(3\,a_2\sspm a_1^2) }. In our case \dstyle{Disc\sspeq -\frac{7^2}{2^2\,3^3} \sspkl 0}, telling that there are three (different) real solutions, in accordance with the explicitly known ones. Therefore, one expects three different constructions for the given right angle zig-zag $B,C,D,E,A$. In the origami version we expect to find three (different) creases $g_1$, $g_2$ and $g_3$  each for folding simultaneously $A$ onto some $A^{\prime}$ on the $x$-axis and $B$ onto some $B^{\prime}$ on the $y$-axis. The three {\it Figures 13,\ 14} and {\it 15} show these solutions.\psn
For all three figures the zeroes of $C(7,x)$ are \dstyle{x\sspeq \frac{\overline{C,\,D}}{R} \sspeq \frac{P_{B_y}}{R} \sspgr 0}. $F\sspeq P_B$ and $G\sspeq P_A$. The folding $A\sspto A^{\prime}$ works like earlier described in connection with {\it Figure 5}. For $B\sspto B^{\prime}$ with $B\,:\, [-2\,R,\, 0]$ (we use the length unit $R$ here) and  $B^{\prime}\,:\, [0,\,\tilde y]$ one  has for the mid-point \dstyle{P_B\,:\, \left[-R,\,\frac{\tilde y}{2}\right]}. With $\gamma\sspeq \angle(D,\,P_B,\,P_A)$, \dstyle{\tan\,\gamma\sspeq \frac{R\sspp \overline{O,\,R_B}}{P_{B_y}}\sspeq \frac{P_{B_y}}{R}}. Hence \dstyle{\overline{O,\,R_B} \sspeq  \frac{P_{B_y}^2}{R}\sspm R\sspeq  \frac{\tilde y^2}{4\,R}\sspm R}. \dstyle{\tan\,\delta\sspeq \tan\left(\frac{\pi}{2} \sspp \gamma\right) \sspeq -\frac{1}{\tan\,\gamma}}. The equation for the crease is \dstyle{g_B\,:\ y\sspeq -\frac{1}{\tan\,\gamma}\,(x\sspm \overline{O,\, R_B }) }.\pn
(Here $x$ is a cartesian varable.) Putting then $\tan\,\beta\sspeq \tan\,\gamma$ (with $\beta\sspeq \angle(A,A^{\prime},O)$) yields \dstyle{\tilde y\sspeq 2\,\frac{a_y}{\hat X} \sspeq \frac{4\,R}{\hat X}} with \dstyle{\hat X\sspdef \frac{X}{R}}, where $X\sspeq R\sspm a^{\prime}_x$. Together with $\overline{0,\,R_A}\sspeq\overline{0,\,R_B} \sspfed r\,R$ one finds  the cubic equation for $\hat X$:\psn
\Beq 
\hat X^3\sspm 4\,\hat X^2\sspm 4\,\hat X \sspp 8 \sspeq 0\, \nonumber
\Eeq
for each of the three figures.\psn
{\bf Figure 13:}\psn
Here $\tilde x\sspeq a^{\prime}_x/R\sspkl 0$, \ie $X\sspeq R\sspp |a_x^{\prime}|\sspgr 0$. \dstyle{x\sspeq \frac{\overline{C,\,P_B}}{R}\sspeq \frac{\tilde y}{2\,R}\sspeq \frac{2}{\hat X}}. The cubic heptagon equation for $x$, given above, is compatible with the cubic equation for $\hat X$. Because the three solutions for $x$ are known from the heptagon (see above), and since here $x\sspgr 1$ one has \dstyle{x\sspeq \rho(7)\sspeq 2\,\cos\left(\frac{\pi}{7}\right) \approx  1.801937736}, corresponding to \dstyle{\hat X\sspeq \frac{2}{x}\sspeq \frac{1}{\cos\left(\frac{\pi}{7}\right)}\sspeq \sqrt{1\sspp \tan\left(\frac{\pi}{7}\right)^2}\approx  1.109916264}.   \psn
{\bf Figure 14:}\psn
$\tilde x\sspeq a^{\prime}_x/R\sspkl 0$, \ie $X\sspeq R\sspp |a_x^{\prime}|\sspgr 0$.
\dstyle{x\sspeq \frac{\overline{C,\,P_B}}{R}\sspeq \frac{\tilde y}{2\,R}\sspeq \frac{2}{\hat X}}. Because $0\sspkl x\sspkl 1$, one has \dstyle{x \sspeq x(7;3)  \sspeq 2\,\cos\left(\frac{3\, \pi}{7}\right) \approx .4450418670},  corresponding to \dstyle{\hat X\sspeq \frac{2}{x}\sspeq \frac{1}{\cos\left(\frac{3\,\pi}{7}\right)}\sspeq \sqrt{1\sspp \tan\left(\frac{3\, \pi}{7}\right)^2}\approx 4.493959217}. \psn
{\bf Figure 15:}\psn
$\tilde x\sspeq a^{\prime}_x/R\sspgr 0$, \ie $X\sspeq R\sspm a_x^{\prime}\sspkl 0$. $\tilde y\sspeq b^{\prime}_y\sspkl 0$. \dstyle{0\sspgr x\sspeq -\frac{\overline{C,\,P_B}}{R}\sspeq \frac{\tilde y}{2\,R}\sspeq \frac{2}{\hat X}}. Hence $x\sspeq x(7;5)\sspeq - 2\,\cos\left(\frac{\pi\,2}{7}\right) \approx -1.246979604$, corresponding to \dstyle{0\sspgr \hat X\sspeq \frac{2}{x}\sspeq -\frac{1}{\cos\left(\frac{2\,\pi}{7}\right)}\sspeq \-\,\sqrt{1\sspp \tan\left(\frac{2\, \pi}{7}\right)^2}\approx -1.603875472}. 
\psn
\parbox{16cm}{
{\includegraphics[height=8cm,width=.5\linewidth]{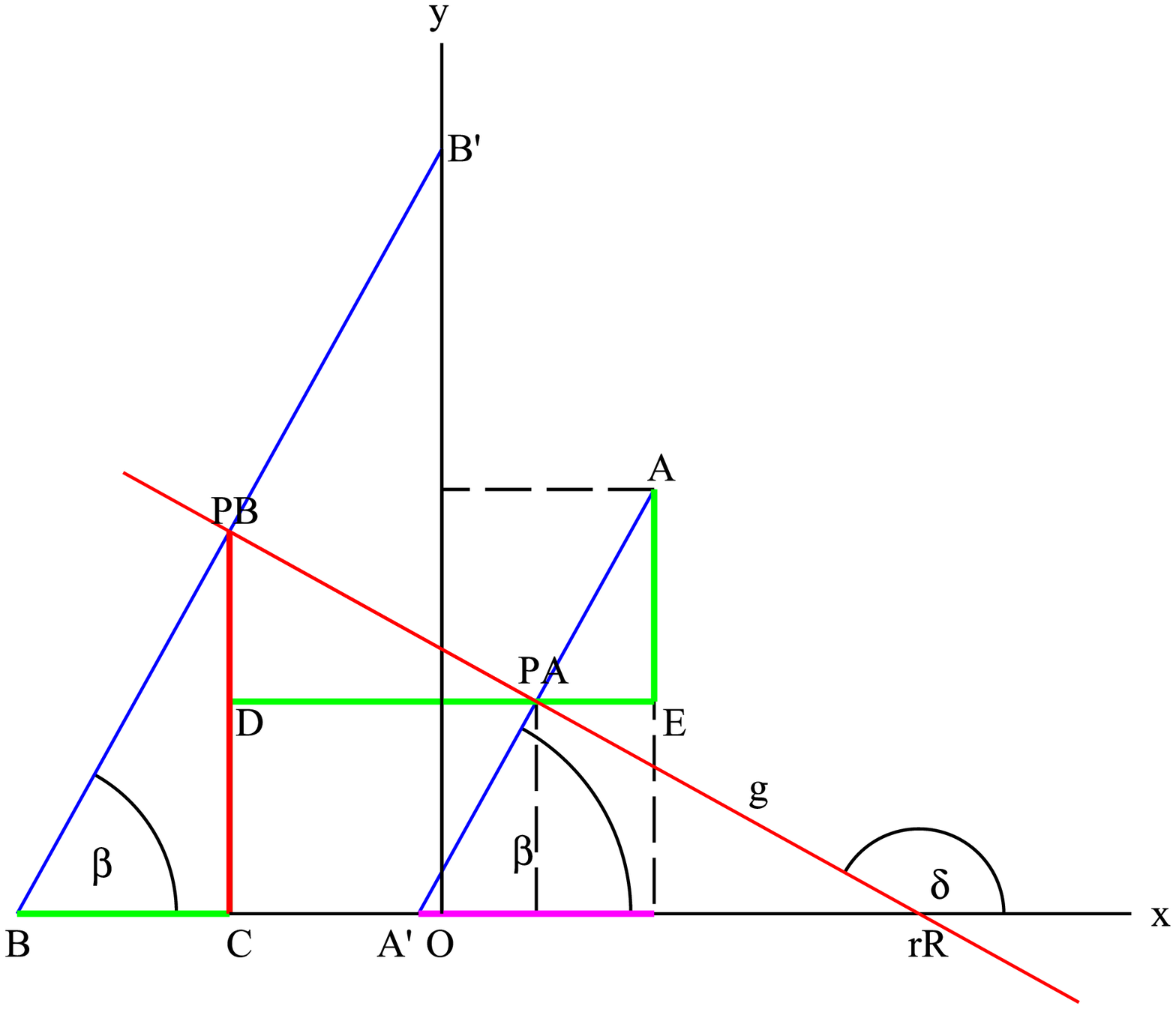}}
{\includegraphics[height=8cm,width=.5\linewidth]{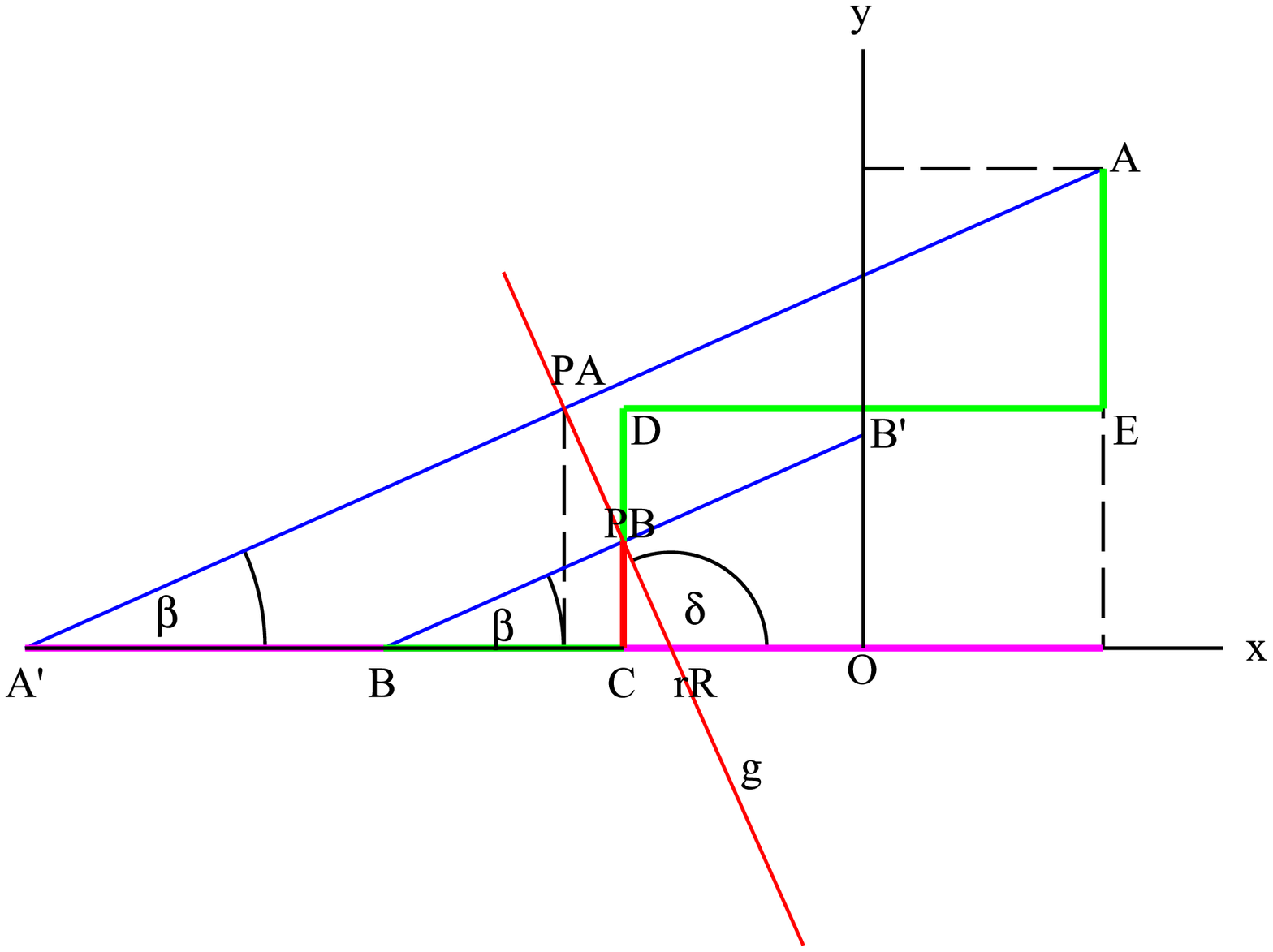}}
}
\psn
\hskip 7cm  Case \dstyle{\alpha\sspeq \frac{\pi}{2}}
\psn
\hskip 1cm Figure 13:\ Heptagon, first origami \hskip 2.5cm   Figure 14:\ Heptagon, second origami \pbn
\pbn
\psn
\parbox{16cm}{\begin{center}
{\includegraphics[height=8cm,width=.5\linewidth]{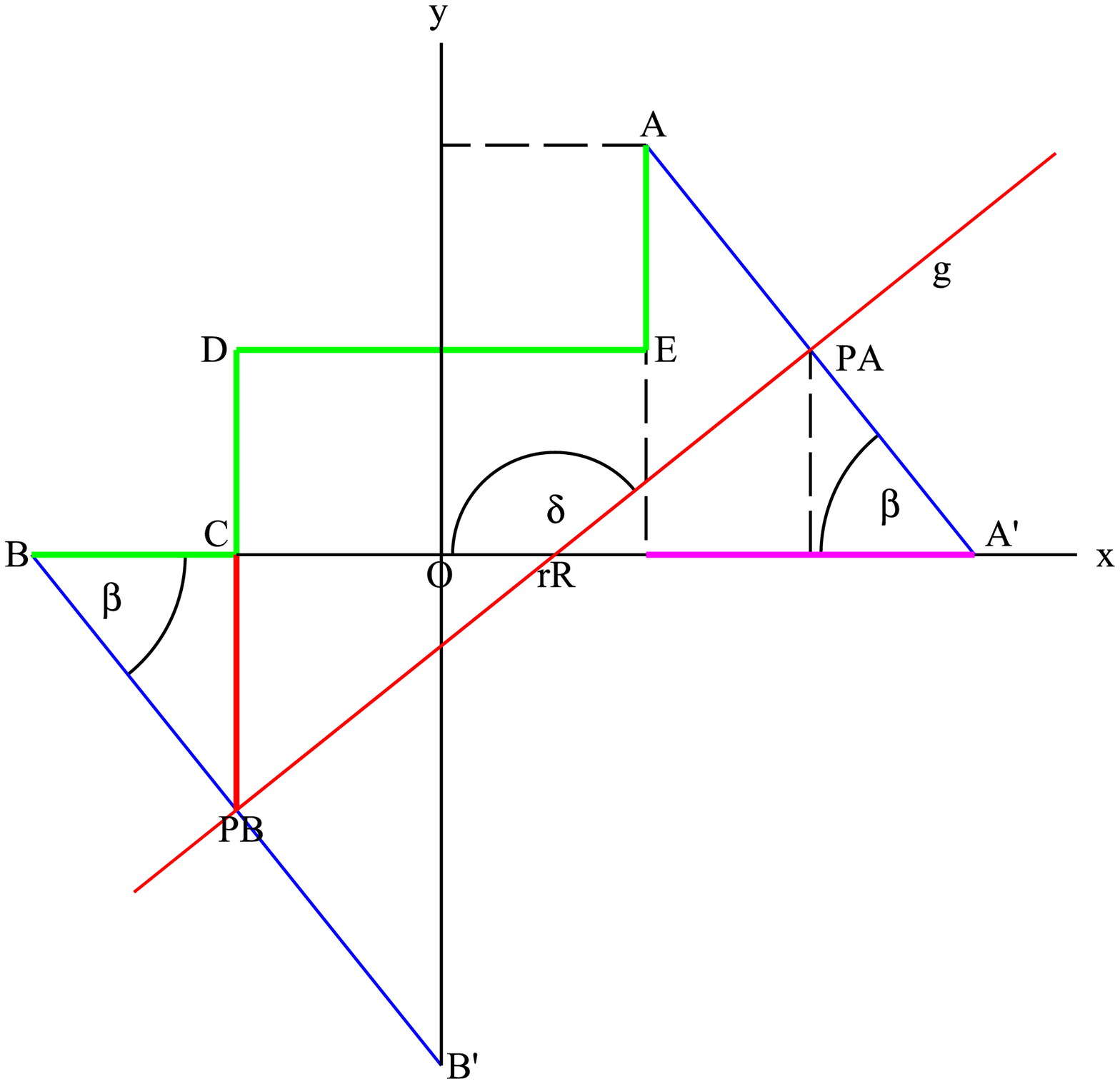}}
\end{center}
}
\psn
\hskip 4cm  \ Figure 15:\ \dstyle{\alpha\sspeq \frac{\pi}{2}},\  heptagon,\ third\  origami 
\pbn
{\bf Application 2: Doubling the cube}\psn
This classical problem cannot be solved by ruler and compass, but with Origami this can be accomplished because one has to solve the third order equation $x^3\sspm 2 \sspeq 0$.  See \cite{Newton}  and {\it Figure 16}. (i)  First one has to find a third of $\overline{A,B}$. This is a standard origami problem solved by finding the intersection point $X$ of the two creases $\overline{A,C}$ and $\overline{A,E}$ where $E$ is found by halving the square, bringing $B\sspto A$ and $C\sspto D$. If the length of the side of the  square $A,\,B,\,C,\,D$ is taken as $1$ (in some length unit) then \dstyle{\overline{B,H}\sspeq \frac{1}{3}\sspeq \overline{C,J}}. 
\pbn
\parbox{16cm}{
{\includegraphics[height=8cm,width=.5\linewidth]{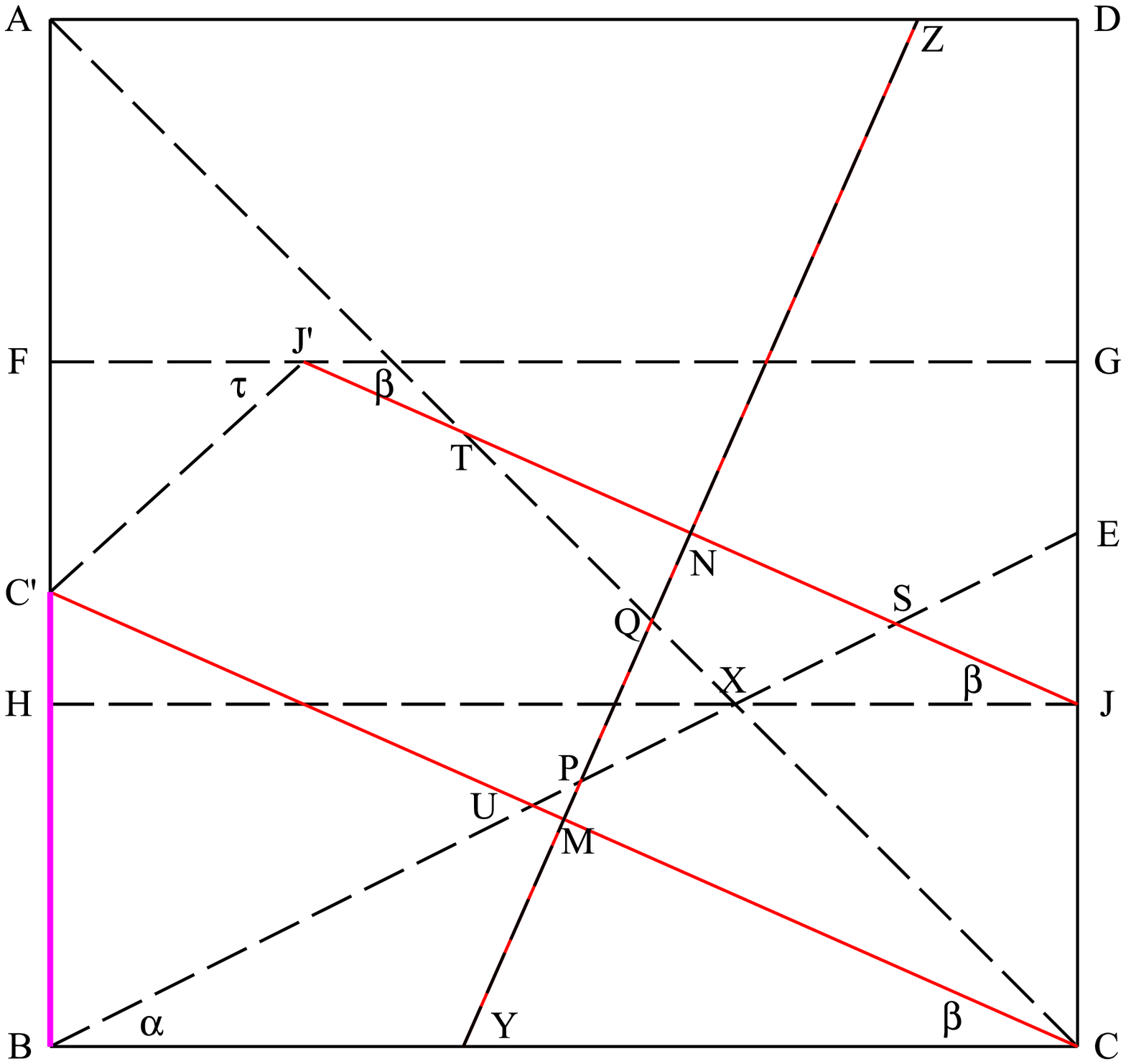}}
{\includegraphics[height=8cm,width=.5\linewidth]{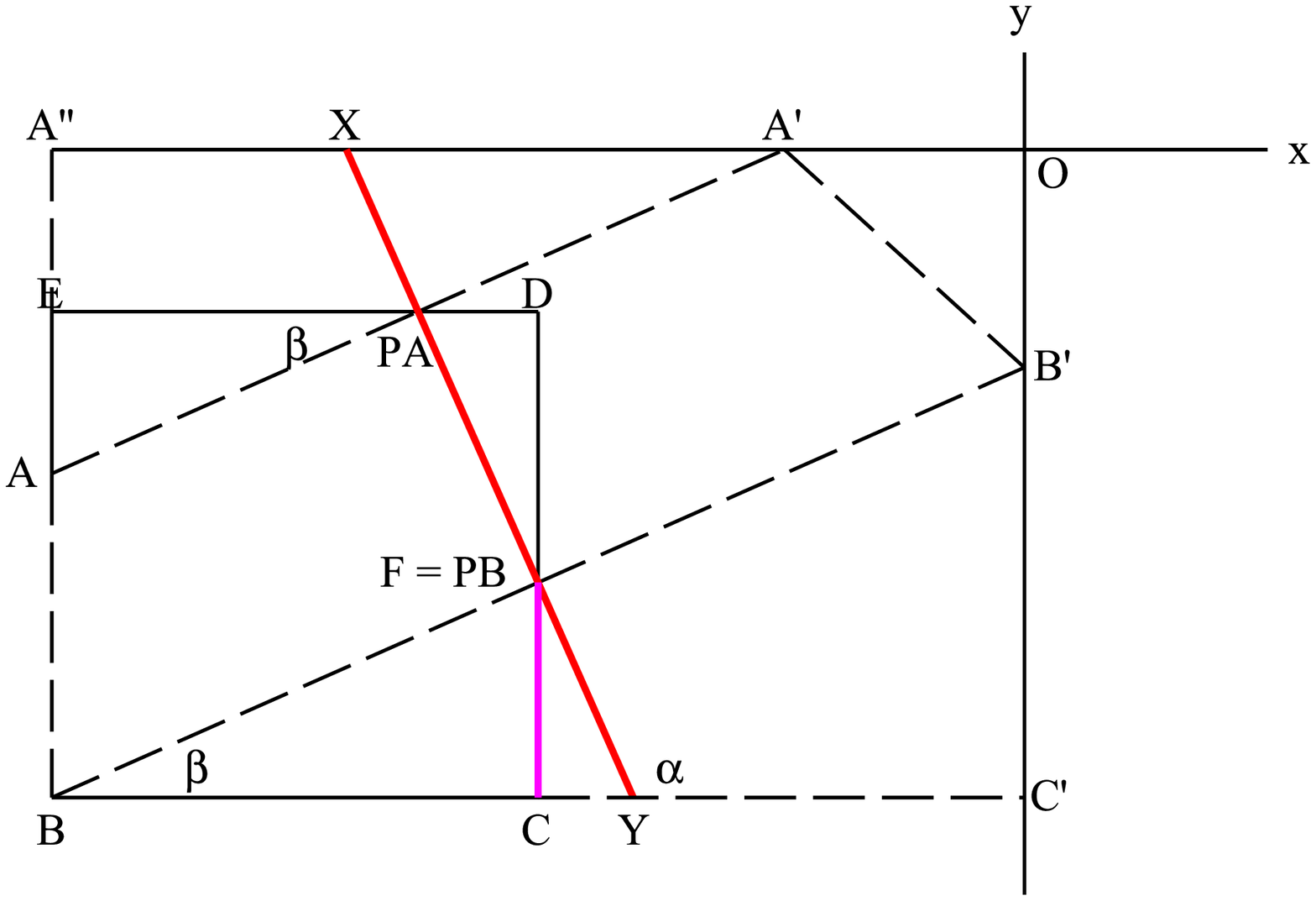}}
}
\psn
 \hskip .1cm Figure 16:\ Doubling the cube; finding $\overline{C^{\prime},B}$ \hskip .5cm Figure 17:\ Doubling the cube; standard version
\pbn
To see this one just has to find the intersection of the two straight lines \dstyle{y\sspeq \frac{1}{2}\, x} and \dstyle{y\sspeq - x +1} finding \dstyle{X:\ \left[\frac{2}{3},\frac{1}{3} \right]}. (ii) Folding $\overline{A,D}$ onto $\overline{H,J}$ will then generate the crease $\overline{F,G}$ which completes the task to divide the square into three equal parts. \dstyle{F:\ \left[0,\frac{2}{3} \right]}, \dstyle{G:\ \left[1,\frac{2}{3} \right]}. (iii) The crucial origami is then to fold at the same time $C\sspto C^{\prime}$ with $C^{\prime}$ on the line $\overline{A,B}$ and $J\sspto J^{\prime}$ with $J^{\prime}$ on the line $\overline{F,G}$. The claim is that \dstyle{x\sspdef\frac{\overline{A,C^{\prime}}}{\overline{C^{\prime},B}}\sspeq 2^{1/3}}, or with \dstyle{s\sspdef \overline{C^{\prime},B}}, \ie \dstyle{x\sspeq \frac{1-s}{s}}, $\bf {s^3\sspm s^2\sspp s \sspm \frac{1}{3}\sspeq 0}$. The discriminant is \dstyle{Disc\sspeq +\frac{1}{3^4}}, showing that there is only one real solution, which is \dstyle{s\sspeq  \frac{1}{3}\,(2^{2/3}\sspm 2^{1/3}\sspp 1)  \sspapprox 0.4424933339}. See \cite{OEIS} \seqnum{A246644} for the decimal expansion of $s$.\pbn
The analytic {\bf  proof} is obtained from looking at the right triangle $\triangle(C^{\prime},\,F,\,J^{\prime}) $ with angle $\tau\sspdef \angle(F\,J^{\prime},\,C^{\prime})$. This angle can be computed by identifying the trapezoid angle $\angle(C^{\prime},\,J^{\prime},\,J)$ with the one \dstyle{\angle(J^{\prime},\,J\, C)\sspeq \beta\sspp \frac{\pi}{2}}. Then \dstyle{\beta \sspp \tau\sspp (\beta\sspp \frac{\pi}{2})\sspeq \pi}, \ie \dstyle{\tau\sspeq \frac{\pi}{2}\sspm 2\,\beta}.  Now \dstyle{\sin\,\tau\sspeq ((1 \sspm s)\sspm \frac{1}{3})/(\frac{1}{3})\sspeq 2\sspm 3\,s}, which is also \dstyle{\sin\left(\frac{\pi}{2}\sspm 2\,\beta\right)\sspeq 2\,\cos^2 \beta\sspm 1}. But \dstyle{\cos\,\beta\sspeq \frac{1}{\sqrt{1\sspp s^2}}} from \dstyle{\tan\,\beta\sspeq s}, and thus \dstyle{\sin\,\tau\sspeq \frac{1-s^2}{1+s^2}}. Equating this with $2\sspm 3\,s$ leads to the claimed equation for $s$, hence the one for the ratio $x$. \hskip 16cm $\square$
\pbn
We list more analytic data for {\it Figure 16} with $s\sspdef \overline{C^{\prime},B},\ b\sspdef \overline{C^{\prime}, M}\sspeq \overline{M,C}$ and  $a\sspdef \overline{J^{\prime}, N}\sspeq \overline{N,J}$ :  \pn
\dstyle{\tan\,\alpha\sspeq \frac{1}{2}},\  $\tan\,\beta\sspeq s/1\sspeq s$,\ \dstyle{\sin\, \beta \sspeq \frac{s}{2\,b}\sspeq \frac{1}{3}/(2\,a)\sspeq \frac{1}{6\,a}}.\ \dstyle{\tan\,\tau\sspeq \frac{1-s^2}{2\,s}} (from the second formula for $\sin\,\tau$ given above).\ \dstyle{J^{\prime}:\, \left[\frac{3\, s-1}{3\,s} ,\frac{2}{3}\right]} from \dstyle{\overline{F,J^{\prime}}\sspeq 1\sspm \overline{J^{\prime},G}\sspeq 1\sspm \frac{1}{3\,\tan\,\beta} \sspeq 1\sspm \frac{1}{3\,s}\sspapprox 0.2466929834}.  \dstyle{M:\, \left[ \frac{1}{2},\,  \frac{s}{2}\right]}, \dstyle{N:\, \left[ \frac{6\,s \sspm 1}{6\,s},\,  \frac{1}{2}\right] \sspeq  \left[ \frac{1\sspp s\sspm s^2}{2},\,  \frac{1}{2}\right]}.\ \dstyle{P:\, \left[\frac{s^2\sspm 1}{s\sspm 2}, \frac{1}{2}\,\frac{s^2\sspm 1}{s\sspm 2}  \right]},\ \dstyle{Q:\, \left[\frac{1\sspp 2\,s\sspm s^2}{2\,(1\sspp s)},\,\frac{1\sspp s^2}{2\,(1\sspp s} \right]}.  \dstyle{S:\, \left[\frac{2}{3}\,\frac{1\sspp 3\,s}{1\sspp 2\,s},\, \frac{1}{3}\,\frac{1\sspp 3\,s}{1\sspp 2\,s}\right]},\ \dstyle{T:\, \left[\frac{1}{3}\,\frac{2\sspm 3\,s}{1\sspm s},\, \frac{1}{3\,(1\sspm s)}\right]}, \pn  
\dstyle{U:\, \left[\frac{2\,s}{1\sspp2\, s},\, \frac{s}{1\sspp2\, s}\right]},\ \dstyle{Y:\, \left[\frac{1\sspm s^2}{2},\, 0\right]}, \ \dstyle{Z:\, \left[\frac{1\sspm s^2}{2}\sspp s,\, 1\right]}. The equation for the  crease $\overline{Y,Z}$ is \dstyle{y\sspeq \frac{1}{s}\,\left(x\sspm \frac{1-s^2}{2} \right)}. \pbn
{\bf Standard version to find $\bf s$ with $\bf s^3\sspm s^2\sspp s \sspm \frac{1}{3}\sspeq 0$} \psn
The cubic equation for $s\sspeq \overline{C^{\prime},B}$ can also be solved geometrically in the standard fashion, similar to finding $x$ in the heptagon case treated above. Here the sign pattern is $+,-,+,-$, which means that the $90^o$ chain pattern is $l,l,l$. This leads to the line chain $B,\,C,\, D,\,E,\, A$ shown in {\it Figure 17} which is identical with $B,\, C,\,D,\, A,\, F$ in {\it Figure 16} (the scale of both figures is different).  In this case the $y$-axis is parallel to $\overline{D,C}$ such that the $x$-coordinates of $B$ becomes $-2$ (for $R\sspeq 1$ length unit), and the $x$-axis is parallel to $\overline{E,D}$ (in {\it Figure 16} this is $\overline{A,D}$) such that the $y$-coordinate of $B$ becomes \dstyle{-\frac{1}{3}\sspm 1\sspeq -\frac{4}{3}}, using \dstyle{\overline{E,A}\sspeq \frac{1}{3}\sspeq \overline{E,A^{\prime\prime}}}. Here the $F$ in the standard geometrical construction (not to be confused with $F$ in {\it Figure 16})  with $\overline{C,F}\sspeq s$  is  shown in {\it Figure 17} .  That is the $[x,y]$ coordinates of this $F=PB$ from the origami construction are \dstyle{\left[-1,-\left(\frac{4}{3}\sspm s\right)\right]}.\ $B\sspto B^{\prime}$ on the $y$-axis with coordinates \dstyle{B^{\prime}:\, [0,\,-\frac{4}{3}\sspp 2\,s]} (from the continuation of the line element  $\overline{B,PB}$ to the $y-$axis). $A\sspto A^{\prime}$ with coordinates  $A^{\prime}$:\ \dstyle{\left[-2\,\left(1\sspm \frac{1}{3\,s}\right) ,\,0\right]} because $\overline{A,A^{\prime}}$  is parallel to $B,B^{\prime}$ with slope $\tan\,\beta\sspeq s$.\ The two  midpoints defining the crease are $F\sspeq P_B$  and $P_A:\, \left[-2\,(1\sspm \frac{1}{6\,s}),\,-\frac{1}{3}\right]$.   The slope of the crease $\overline{X,Y}$, shown in {\it Figure 17}, is \dstyle{\tan\,\alpha\sspeq -\frac{1}{\tan\, \beta}\sspeq-\frac{1}{s}}. The equation for the crease is \dstyle{y\sspeq -\frac{1}{s}\,\left(x\sspp \frac{s^2\sspp 6\,s\sspm 1}{3\,s}\right)}. The data for the trapezoid $A,\,B,\,B^{\prime},\, A^{\prime}$ is: \dstyle{a\sspdef \overline{A,PA}\sspeq \overline{PA,A^{\prime}} \sspeq \frac{\sqrt{1+s^2}} {3\,s}\sspapprox 0.8237612353}, and \dstyle{b\sspdef \overline{B,PB}\sspeq \overline{PB,B^{\prime}} \sspeq \sqrt{1+s^2}\sspapprox 1.093526566}.\  \dstyle{\overline{A,B}\sspeq \overline{A,A^{\prime\prime}}\sspeq \frac{2}{3} }.
\pbn
{\bf Completing the task of doubling a given cube}\psn
Up to now we have only found the doubling of the cube with side length $s$. In {\it Figure 16} we had $2^{1/3}\, \overline{C^{\prime},B}\sspeq \overline{C^{\prime},A}$, \ie $2\,s^3\sspeq (1\sspm s)^3$. In {\it Figure 17} we had $2^{1/3}\,\overline{PB,C}\sspeq \overline{D,PB}$.  For the decimal expansion of $2^{1/3}$ see \cite{OEIS} \seqnum{A002580}. But the task is to double a cube with given side length $L$. If one takes as length of the side of the square \dstyle{R\sspeq \frac{L}{s}}, then $L\sspeq \overline{B,C^{\prime}}$ and  $M\sspeq \overline{C^{\prime},A}\sspeq 2^{1/3}\,L$ is the side length for the doubled cube. However, we first have to find \via\, origami $1/s\sspapprox 2.259921051$. But this can be achieved by considering  the parallel to $\overline{C,C^{\prime}}$ through $A$. This parallel will hit the continuation of $\overline{B,C}$ on some point $C^{\prime\prime}$ with coordinates \dstyle{\left[\frac{1}{s},\ 0\right]}  (origin at $B$, $x$-axis along $\overline{B,C}$ and $y$-axis along $\overline{B,A}$). See {\sl Figure 18}. This means that if we take the length scale $R\sspeq L$ the searched length $M\sspeq 2^{1/3}\,L$ for the doubled cube is given by $\overline{C,C^{\prime\prime}}$ which is \dstyle{\left(1-\frac{1}{s}\right)\,L}.  It is easy to find the parallel $g_2$ in {\it Figure 18} by origami. First find $g1$, the crease perpendicular to crease $g$ through the point $A$ (this can be done, as explained in the introductory remarks; axiom 4 or IV). Then find $g2$ as the crease perpendicular to $g1$ through point $A$. Finally the square $C,C^{\prime\prime},D^{\prime\prime},A^{\prime\prime}$ can be completed.\psn
The coordinates of some points are:   $C:\ [L, 0]$,\  $C^{\prime}:\ [0,s\,L]$,\ \dstyle{C^{\prime\prime}:\ \left[\frac{L}{s}\sspeq (1\sspp 2^{1/3})\, L),\, 0\right]},\  $V:\ [-s\,L,0]$,\  \dstyle{W:\ \left[- \frac{s\,(1-s)}{1\sspp s^2},\,L\, \frac{s\,(1+s)}{1\sspp s^2}\,L\right]}, \ \dstyle{D^{\prime\prime}:\ \left[\frac{L}{s},\, 2^{1/3}\,L  \right]}.
\psn
\parbox{16cm}{\begin{center}
{\includegraphics[height=8cm,width=.5\linewidth]{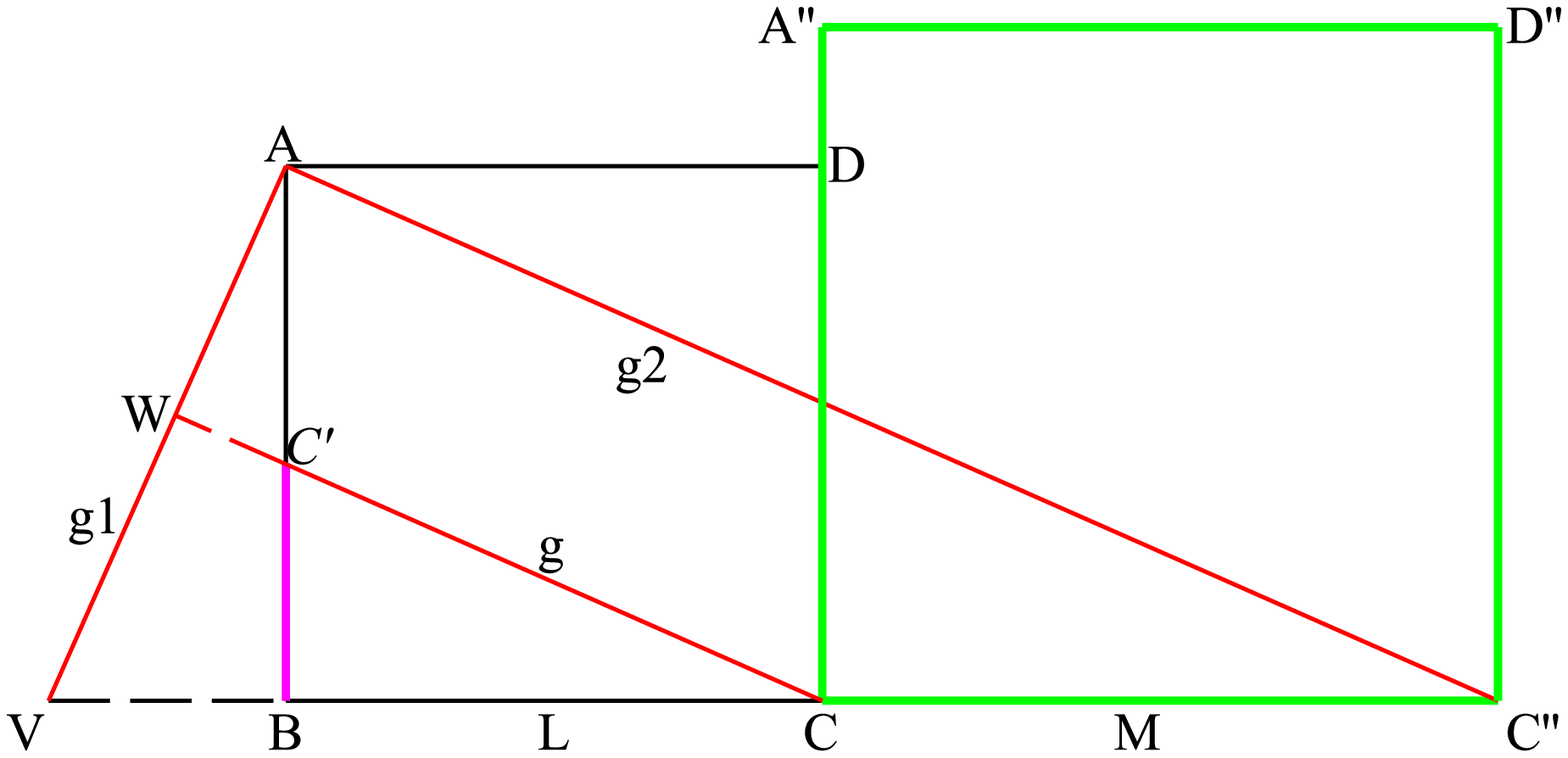}}
\end{center}
}
\psn
\hskip 4cm  \ Figure 18: Doubling the cube: finding $M \sspeq 2^{1/3}\,L$
\pbn
{\bf Application 3: Trisection of an angle}\psn
This is another classical problem not solvable with ruler and  compass but with origami. See \cite{Husini}, \cite{Huzita2}, \cite{Newton}. We first discuss the origami shown in \cite{Huzita2} {\it Fig. 1}, pp. 204-5, and also in \cite{Newton}.  See the present {\it Figure 19}, where the angle is $\alpha\sspeq \angle (P,B,C)$. Because the origami solution will be based on a cubic equation with three real roots, the question of the meaning  of the other two roots arises. The answer can be found in \cite{Huzita2}: the origami prescription for trisecting a given angle is not unique and the other two solutions correspond to trisecting the angle $\pi\sspm \alpha$ and $\pi\sspp \alpha$. This will be treated at the end. of this section.\psn
In the square $(A,\,B,\,C,\,D)$ the point $P$ on $\overline{A,D}$ defines the angle $\alpha \sspeq \angle {PBC}$ to be trisected and the crease $g1$. Then an arbitrary horizontal crease $g2$ defining points $E$ and $F$ with distance $2\,h \sspkl 1$ from the base line $\overline{B,C}$ is folded. The dashed crease $g3$ bringing $B$ onto $E$  and $C$ onto $F$ has then a distance $h$ from the base line. The crucial folding $g$ is then to bring point $E$ onto $E^{\prime}$ on crease $g1$ and simultaneously point $B$ onto $B^{\prime}$ on crease $g3$. This will also bring point $G$ to  $G^{\prime}$.  The continuation of $\overline{B,B^{\prime}}$ will intersect the line $\overline{D,C}$ at a point $V$, defining crease $g4$ .The intersection point of crease $g$ with crease $g3$ is called $Y$. This defines the blue crease $g5$ with line segment $\overline{B,Y}$ crossing the continuation of the top line $\overline{A,D}$ at a point $Q$ (depending on the choice of $P$ and $h$ this point $Q$ could also lie on $\overline{A,D}$, \eg for $\alpha\sspeq 70^o$ and $h\sspeq 0.2$). The claim is now that the blue crease $g5$ and the crease $g4$ trisect the angle \dstyle{\alpha$ with $\sigma\sspeq \frac{\alpha}{3}}. One also shows that  $G^{\prime}$ lies on the blue crease $g5$.\psn
The following  analytic data is given for a coordinate system with origin at $B$, the $x-$axis along $\overline{B,C}$ and the $y$-axis along $\overline{B,A}$. The input quantities are  $R$, the length of the square in some length unit, \dstyle{\overline{A,P}\sspeq x\sspeq \frac{1}{tan\,\alpha}}, with input $\alpha$ in radians, and $\overline{B,\,G}\sspeq h\, R$ :\,\pn
$\beta\sspeq \frac{\pi}{2}\sspm \sigma$,\ $\overline{B,\,C} \sspeq R$,\  $E:\, [0,2\,h\,\,R]$,\ $G:\, [0,h\,\,R]$, \ \ $F:\, [R,\,2\,h\,R]$,\ $H:\, [R,h\,\,R]$,\ \dstyle{B^{\prime}:\, \left[\frac{h}{\tan\,\sigma}\,R,\, h\,R \right] },\pn
\dstyle{E^{\prime}:\, \left[h\,\frac{\cos\,\alpha}{\sin\,\sigma}\,R,\, h\,\frac{\sin\,\alpha}{\sin\, \sigma}\,\,R \right] },\ \ \dstyle{G^{\prime}:\, \left[\frac{1\sspp \cos\,(2\,\sigma)}{\tan\,(2\,\sigma)}\,h\,R, \,(1\sspp \cos(2\,\sigma))\,h\,R \right]},\ \ \dstyle{Y:\, \left[\frac{h}{\tan(2\,\sigma)}\,R,\,h\, R \right]},\pn
 \dstyle{X:\, \left[\frac{h}{2\,\tan\, \sigma}\,R, \frac{h}{2}\, R \right]},\ \ \dstyle{L:\, \left[0,\,\frac{h}{\sin(2\,\sigma)\,tan\, \sigma}\, R \right]},\ \  \dstyle{J:\, \left[\frac{h}{\sin(2\,\sigma)}\,R,\,0\right]}, \pn
\dstyle{K:\, \left[\sin(2\,\sigma)\,\frac{(L_y/R \sspm h)}{2}\,R,\,(h\sspp \sin^2\,\sigma\,(L_y/R\sspm h)\, R \right]},\ \pn
\dstyle{Z:\, \left[\sin(2\,\sigma)\,\frac{(L_y/R \sspm 2\,h)}{2}\,R,\,(2\,h\sspp \sin^2\,\sigma\,(L_y/R\sspm 2\,h)\, R \right]},\ \ \dstyle{Q:\, \left[\frac{1}{\tan(2\,\sigma)}\,R,\,R\right]},\ \ \dstyle{V:\, [R,\, R\, \tan\sigma]}. \pn
The equations for the creases are ($x$ is here the abscissa): $g1:\, y\sspeq \tan(\alpha)\, x$,\ \ $g2:\ y\sspeq 2\,h$,\ \ $g3:\ y\sspeq h$,\ \ \dstyle{g:\ y\sspeq \tan\left(\frac{\pi}{2}\sspp \sigma\right)\,\left(x\sspm \frac{h}{\sin(2\,\sigma)}\right)},\ \ $g4:\ y\sspeq \tan(\sigma)\, x$,\ \ $g5:\ y\sspeq \tan(2\,\sigma)\, x$.\pn
Some lengths in the trapezoid $(B,\, B^{\prime},\,E^{\prime},\,E)$ are: $\overline{E^{\prime},G^{\prime}}\sspeq \overline{E,G}\sspeq h\,R\sspeq \overline{G^{\prime},B^{\prime}}\sspeq \overline{G,B}$,\  \dstyle{b\sspdef \overline{B,X}\sspeq \overline{X,B^{\prime}}\sspeq \frac{h}{2\,\sin\,\sigma}\, R},\ \ \dstyle{e\sspdef \overline{E,Z}\sspeq \overline{Z,E^{\prime}}\sspeq \frac{\sin\, \alpha\sspm 2\,sin\,\sigma}{2\,\sin^2 \sigma}\,h\, R},\ \ \dstyle{\overline{B,Y} \sspeq \frac{h}{\sin\,(2\,\sigma)}\,R}. \pbn
Now to the {\bf proof} of the trisection of $\alpha \sspeq \angle (P,B,C)$. Name the three angles, called $\sigma$ in {\it Figures 19} and {\it 20}, as follows. $\tau\sspdef \angle (X,B,J)$, $\sigma\sspdef \angle (Y,B,X,)$ and $\eta\sspdef \angle (E^{\prime},B,G^{\prime})\sspeq \alpha\sspm (\tau\sspp \sigma)$. We want to show that $\tau\sspeq \sigma$ and $\eta\sspeq \sigma$ which implies $\alpha\sspeq 3\,\sigma$. \pn
Consider the  angle $\varepsilon\sspdef \angle(G^{\prime},Y,B^{\prime})$. (Note that at this stage it is not yet clear that  $\varepsilon \sspeq \angle (Y,B,J)$ which would immediately show that $\sigma \sspeq \tau$.  This is because it is not yet clear that the line $\overline{Y,G^{\prime}}$ (obtained from folding along  $g$ where $Y$ is the intersection of $g$ with $\overline{G,H}$) really continues to point $B$. For this one has to prove $\sigma\sspeq \tau$.) Because $\angle(G,Y,L)\sspeq \angle (B,J,L) \sspfed \beta\sspeq \frac{\pi}{2}\sspm \tau$, and also $\angle(K,Y,G^{\prime})\sspeq \angle(K,Y,G) \sspeq \beta$ from the folding along $g$, we have $\pi\sspeq 2\,\beta\sspp \varepsilon$, or $\varepsilon\sspeq \pi\sspm 2\,\beta\sspeq 2\,\tau$. The proof that $\sigma \sspeq \tau$ is done by starting with $\sigma\sspeq \angle(X,B^{\prime},Y)$ from the folding along $g$. Also \dstyle{\angle(Y,B^{\prime},G^{\prime})\sspeq \frac{\pi}{2}\sspm \varepsilon} because of the folding along $g$ the right angle $\angle (Y,G,E)$ appears also as $\angle (Y,G^{\prime},E^{\prime})$ and $E, G^{\prime}$ and $B^{\prime}$ are on a straight line, like $E, G$ and $B$.  Because \dstyle{\angle(X,B^{\prime},Y)\sspeq \frac{\pi}{2}\sspm \tau} we have from the right angle $\angle(J,B^{\prime},E^{\prime})$  also $\angle (J,B^{\prime},G^{\prime})\sspeq \frac{\pi}{2}\sspeq \tau \sspp \sigma \sspp (\frac{\pi}{2}\sspm \varepsilon)$ \ie  $0\sspeq \tau\sspp \sigma-2\,\tau$ or $\sigma \sspeq \tau $.  This proves that the point  $G^{\prime}$ lies on the straight line connecting  $B$ and  $Y$, defining the crease $g5$. Finally, $\eta\sspeq \sigma$ because the blue line $\overline{B,G^{\prime}}$ is the height, let its length be $k$, in $\triangle(B,B^{\prime},E^{\prime})$ and $tan\,\sigma\sspeq \frac{h}{k}\sspeq tan\,\eta$ (because $\overline{G^{\prime},B^{\prime}}\sspeq  \overline{G,B}\sspeq \sspeq \overline{G^{\prime},E^{\prime}}\sspeq h$ from folding along $g$). This implies that this triangle is isosceles, \ie $\overline{B,E^{\prime}}\sspeq \overline{B,B^{\prime}}$. \hskip 16cm $\square$\pbn
\pbn
\parbox{16cm}{
{\includegraphics[height=8cm,width=.5\linewidth]{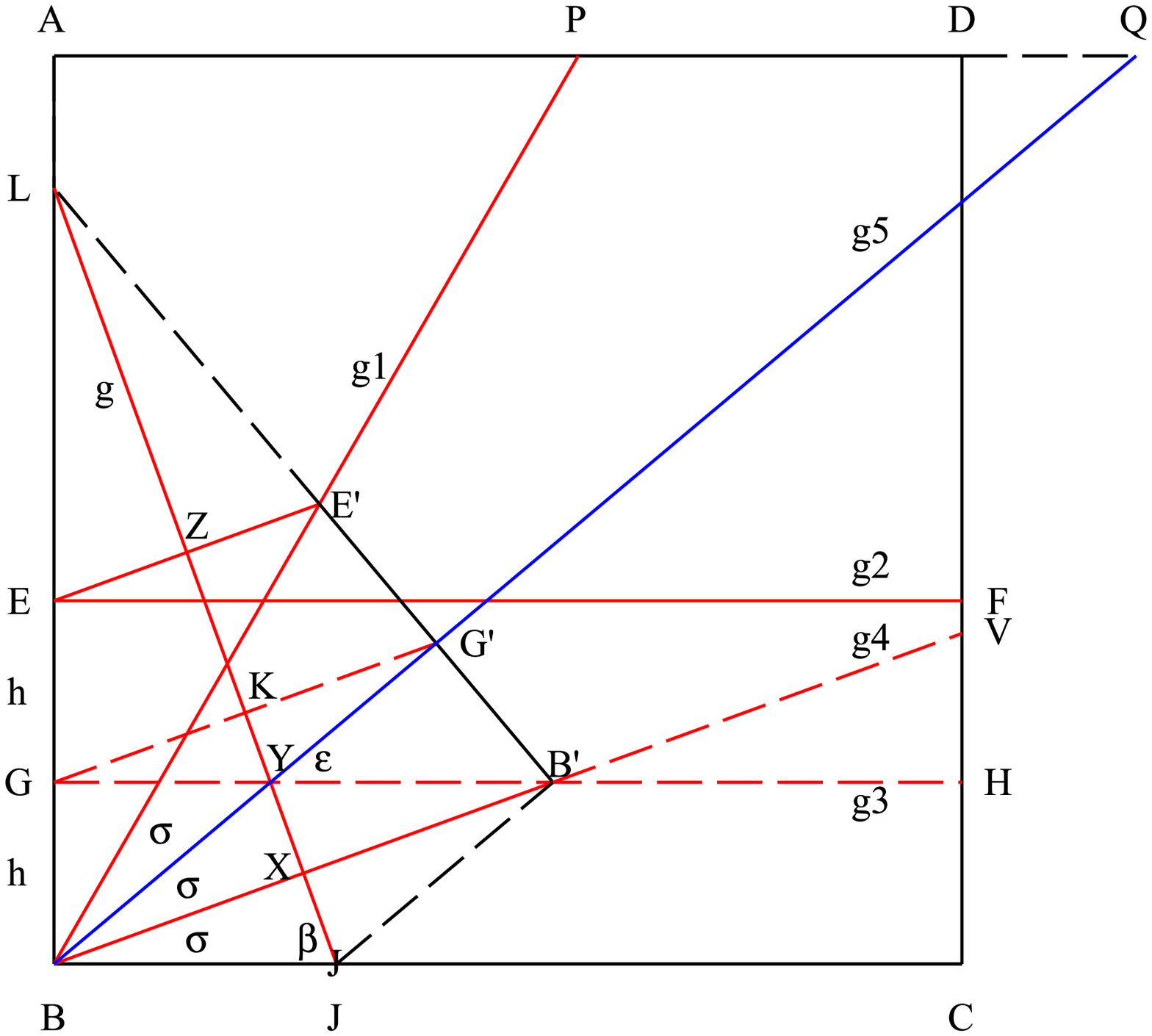}}
{\includegraphics[height=8cm,width=.5\linewidth]{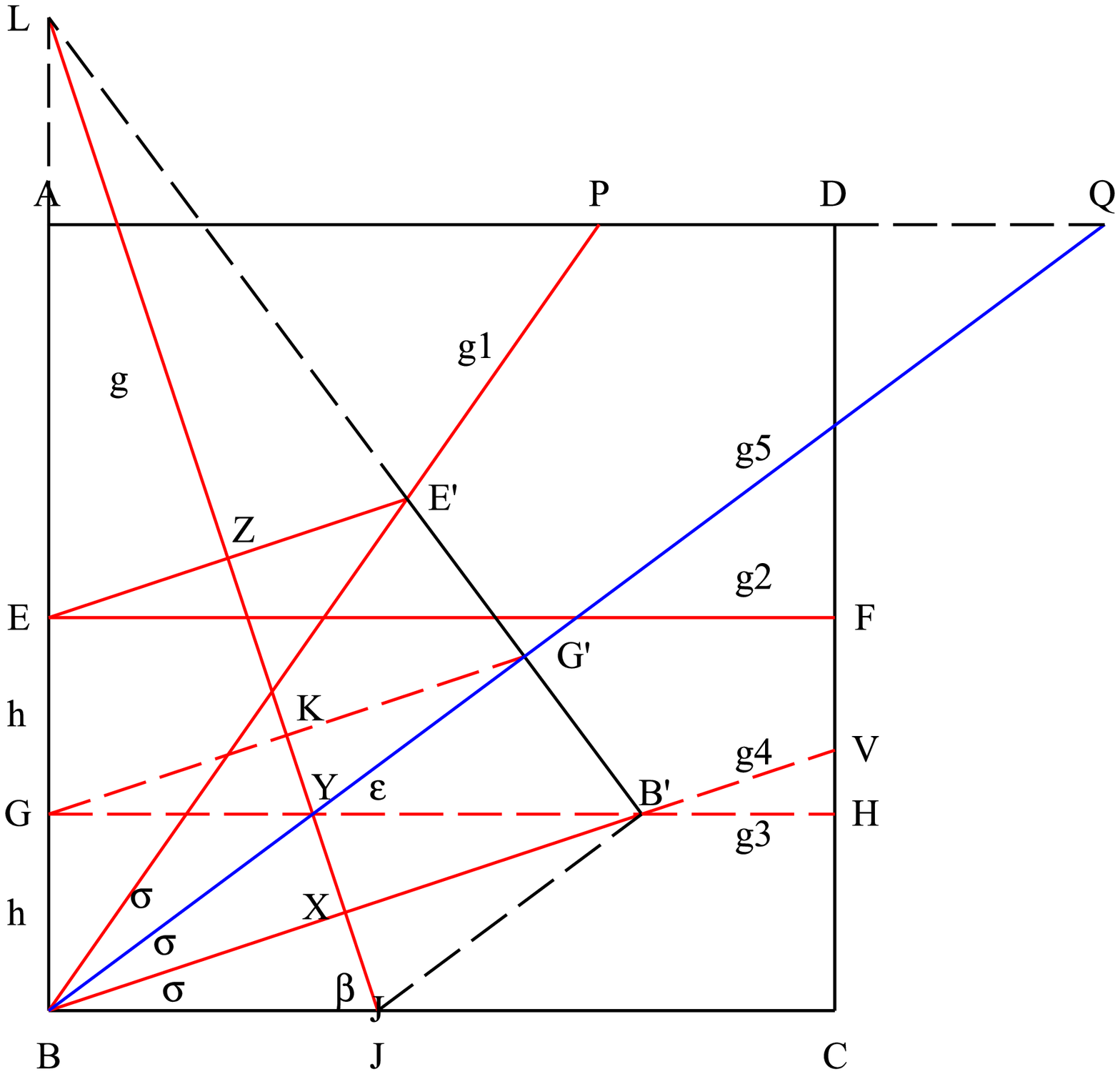}}
}
\psn
 \hskip .1cm Figure 19:\ Trisecting an angle \hskip 3.5cm Figure 20:\ Trisecting an angle\pn
  \hskip 1.5cm $\alpha \sspeq 60^o$ and $h\sspeq 0.2$\hskip 4.5cm  $\alpha\sspeq 55^o$ and $h\sspeq 0.25$\psn
\pbn
As mentioned above, in \cite{Huzita2} the cubic equation relevant for this origami trisection of an angle $\alpha$ is identified and the question of the other two roots is answered.   
The crucial prescription to fold $B \sspto B^{\prime}$ and, at the same time $E\sspto E^{\prime}$, has three solutions, corresponding to the three real solutions of the cubic equation governing this folding. In {\it Figure 21} the origami for the trisection of the angle $\alpha$ is repeated but now  $X1 \sspdef \overline{G, B^{\prime}_x}$ has been shown as a fat line segment (in blue). If one uses the results from above for this case one will find that the dimensionless \dstyle{\hat X1\sspdef  \frac{X1}{R}\sspeq \frac{\overline{G,B^{\prime}}}{R} \frac{h}{tan\, \frac{\alpha}{3}}}  satisfies the following cubic equation.
\Beq
{\hat X}^3 - \frac{3\,h}{\tan\,\alpha }\, {\hat X}^2 \sspm  3\,h^2\,{\hat X} \sspp \frac{h^3}{\tan\,\alpha}\sspeq 0\ .  \nonumber
\Eeq
The discriminant of this cubic is \dstyle{Disc \sspeq -\frac{h^6\,(1\sspp \tan(\alpha)^2)^2}{\tan(\alpha)^4}}, hence negative, therefore there are three different real solutions. The other two solutions $X2$ and $X3$ are shown in the {\it Figures 22} and {\it 23}, respectively.
 \psn
\parbox{16cm}{\begin{center}
{\includegraphics[height=8cm,width=.5\linewidth]{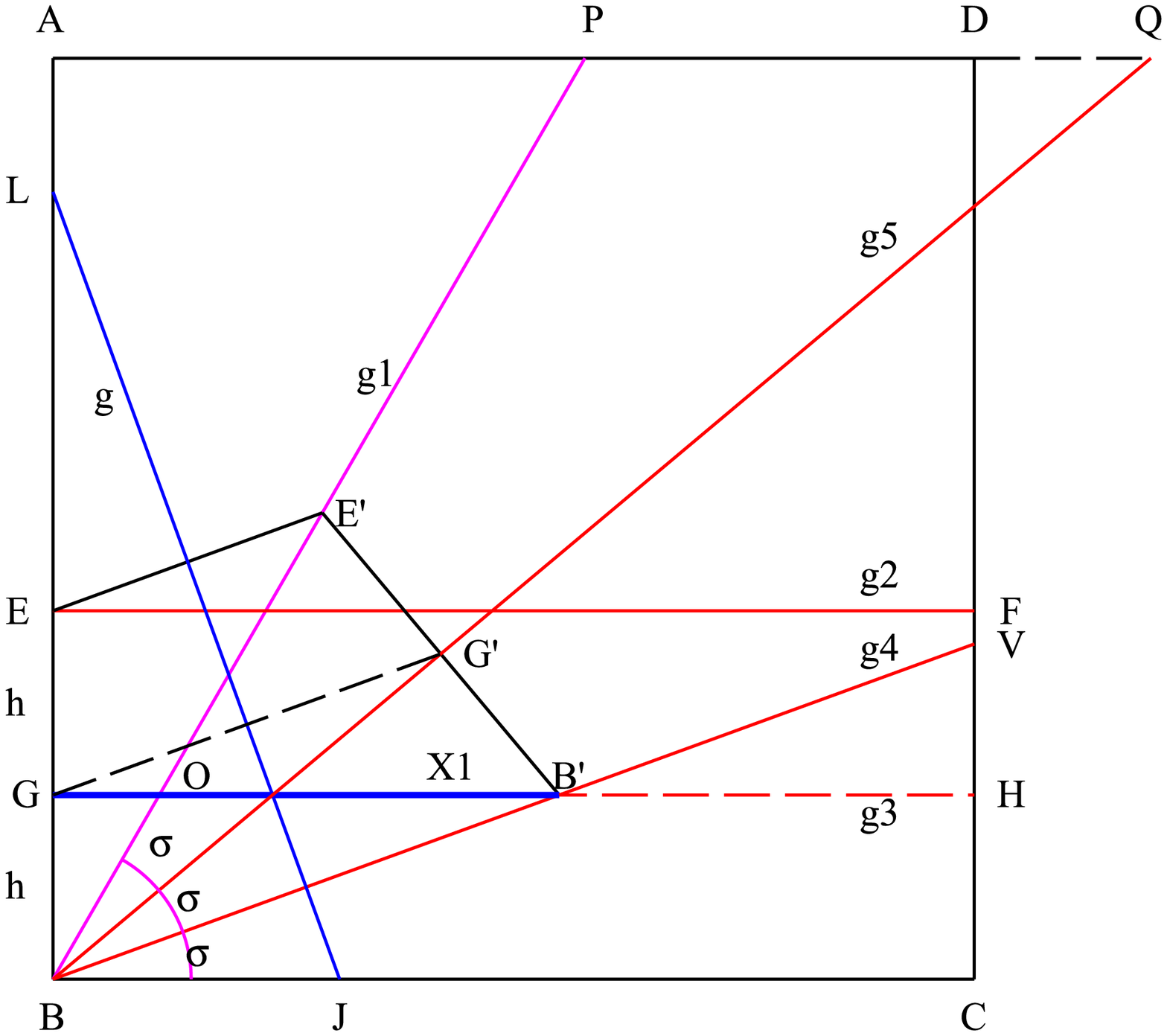}}
\end{center}}
\psn
 \hskip 4cm Figure 21:\ Solution $X1$: Trisecting an angle $\alpha$\pn
\hskip 7cm  $\alpha\sspeq 60^o$ and $h\sspeq .2$
\psn
\parbox{16cm}{
{\includegraphics[height=8cm,width=.5\linewidth]{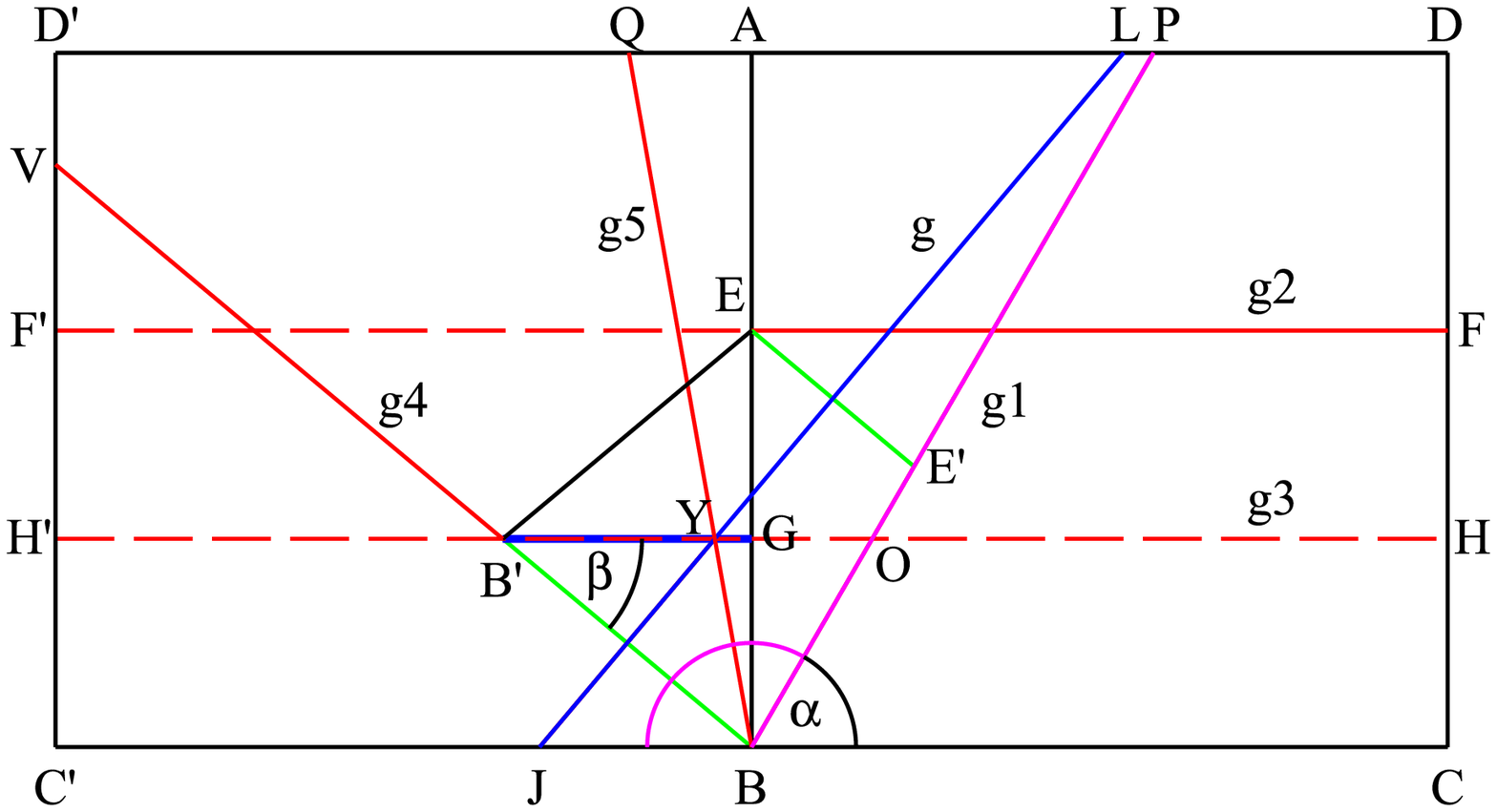}}
{\includegraphics[height=8cm,width=.5\linewidth]{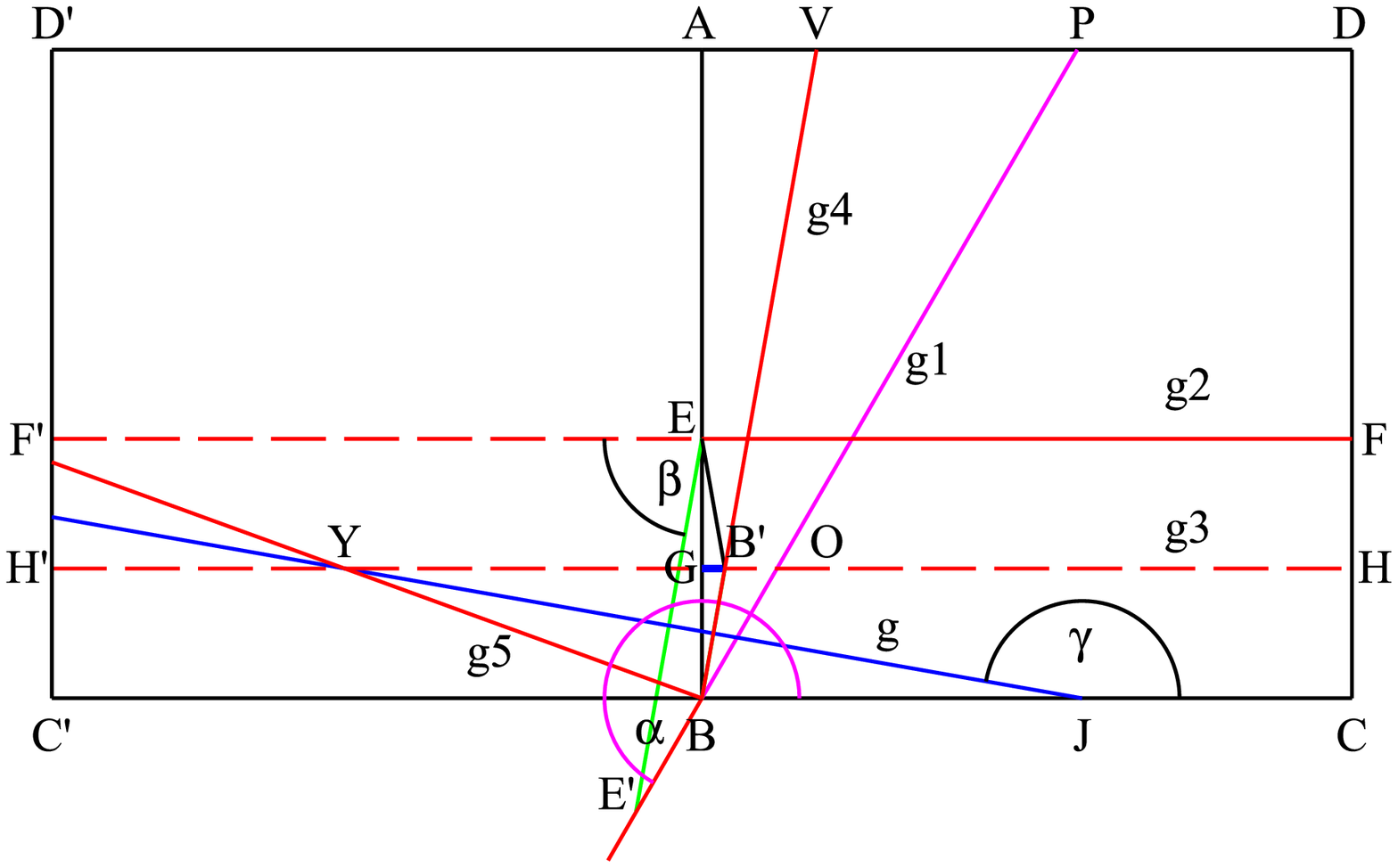}}
}
\psn
 \hskip .6cm Figure 22:\ $X2$,\ Trisecting an angle $\pi \sspm \alpha$ \hskip 2cm Trisecting an angle $\pi\sspp \alpha$ \pn
\hskip 1.9cm  $180^0 \sspm 60^o \sspeq 120^o$, $h\sspeq .3$ \hskip 2.9cm $180^0 \sspp 60^o\sspeq 240^o$, $h\sspeq .2$
\pbn
Finally, we list some analytic data for the three trisection {\it Figures}. The data for {\it Figure 21} has been given already.  The origin is taken at $B$ with the $x-$axis along $\overline{B,C}$ and the $y$-axis along $\overline{B,A}$. Note that in the origami for the general cubic equation treated above the $y^{\prime}$ axis was along the crease $g1$ and the origin was there at \dstyle{O:\ \left[ \frac{h\,R}{tan\,\alpha},\,h\,R\right]}.\ \dstyle{\sigma\sspeq \frac{\alpha}{3}}, as used also above .\psn
{\bf Figure 21:\ \dstyle{{\bf \hat X1 \sspeq \frac{X1}{R}\sspeq \frac{\overline{G,B^{\prime}}}{R}\ (\approx 0.5495\ for\  R\sspeq 1, {\boldsymbol{\alpha}} \sspeq 60^o, h\sspeq .2)}} }\psn
\Beq
 0\sspkl \frac{h}{\hat X1}\sspeq tan\,\sigma\sspeq \frac{e_y^{\prime}\sspm 2\,h}{e_x^{\prime}} \sspeq \frac{\sin\,\alpha \sspm 2\,\sin\, \sigma}{\cos\, \alpha}\ . \nonumber
\Eeq
In the first equation the mapping $B\sspto B^{\prime}$ and in the second one $E\sspto E^{\prime}$ has been considered. Thus \dstyle{\hat{X1}\sspeq \frac{h}{\tan\,\sigma}}, which is in the {Figure 21} shown for $R\sspeq 1$, $\alpha\sspeq 60^o$ and $h\sspeq .2$ which has the above given value.\psn
{\bf Figure 22:\ \dstyle{{\bf \hat X2 \sspeq \frac{X2}{R}\sspeq -\frac{\overline{G,B^{\prime}}}{R}\, (\approx -.23835\ for\  R\sspeq 1,{\boldsymbol  {\alpha}}\sspeq 60^o, h\sspeq .3)}} }\psn
\dstyle{\beta \sspeq \frac{\pi\sspm \alpha}{3}},\ $\overline{B,G}\sspeq h\,R\sspeq \overline{G,E}$, \dstyle{\tan\, \beta\sspeq \tan\left(\frac{\pi-\alpha}{3}\right)\sspeq \frac{h}{|\hat X2|} \sspeq \frac{2\,h\,R\sspm e^{\prime}_y}{e^{\prime}_x}},\pn  
\dstyle{\overline{B,B^{\prime}} \sspfed 2\,b\,R\sspeq \frac{h\,R}{\sin\, \beta}} , \dstyle{\overline{E,E^{\prime}}\sspfed 2\,e\, R \sspeq \frac{e^{\prime}_x}{\cos\,\beta}}, \dstyle{B^{\prime}:\, \left[-\frac{h\,R}{\tan\,\beta},\, h\,R\right]},\pn
\dstyle{E^{\prime}:\left[2\,h\,R\,\frac{\cos\,\alpha\,\cos\, \beta}{\sin\,(\alpha\sspp \beta)},\, 2\,h\,R\,\frac{\sin\,\alpha\,\cos\, \beta}{\sin\,(\alpha\sspp \beta)}  \right]}, \dstyle{J:\ \left[-\frac{h\,R}{\sin(2\,\beta)},\, 0\right]}.\pn
The equations for the creases are: \dstyle{g:\ y\sspeq \frac{1}{\tan\,\beta},(x\sspp \frac{h\,R}{\sin(2\, \beta)})},\ \dstyle{g1:\ y \sspeq (\tan\,\alpha)\, x},\pn
 \dstyle{g4:\ y\sspeq (\tan\,\beta)\,x,\  g5:\ y\sspeq -\tan(2\, \beta)\, x }.\pbn
{\bf Figure 23:\ \dstyle{{\bf \hat X3 \sspeq \frac{X3}{R}\sspeq \frac{\overline{G,B^{\prime}}}{R}\, (\approx +0.0353\ for\  R\sspeq 1, {\boldsymbol {\alpha}}\sspeq 60^o, h\sspeq .2)}} }\psn
\dstyle{\beta \sspeq \frac{\pi\sspp \alpha}{3}},\ \dstyle{\gamma\sspeq \frac{\pi}{2}\sspp \beta},\  $\overline{B,G}\sspeq h\,R\sspeq \overline{G,E}$, \dstyle{\tan\, \beta\sspeq \tan\left(\frac{\pi+\alpha}{3}\right)\sspeq \frac{h}{\hat X3} \sspeq \frac{2\,h\,R\sspp |e^{\prime}_y|}{|e^{\prime}_x|}},\pn  
\dstyle{\overline{B,B^{\prime}} \sspfed 2\,b\,R\sspeq \frac{h\,R}{\sin\, \beta}} , \dstyle{\overline{E,E^{\prime}}\sspfed 2\,e\, R \sspeq \frac{|e^{\prime}_x|}{\cos\,\beta}}, \dstyle{B^{\prime}:\, \left[\frac{h\,R}{\tan\,\beta},\, h\,R\right]},\pn
\dstyle{E^{\prime}:\left[2\,h\,R\,\frac{\cos\,\alpha\,\cos\, \beta}{\sin\,(\beta\sspm \alpha)},\, 2\,h\,R\,\frac{\sin\,\alpha\,\cos\, \beta}{\sin\,(\beta\sspp \alpha)}  \right]}, \dstyle{J:\ \left[\frac{h\,R}{\sin(2\,\beta)},\, 0\right]}.\pn
The equations for the creases are: \dstyle{g:\ y\sspeq -\frac{1}{\tan\,\beta}\,\left(x\sspp \frac{h\,R}{\sin(2\, \beta)}\right)},\ \dstyle{g1:\ y \sspeq (\tan\,\alpha)\, x},\pn
 \dstyle{g4:\ y\sspeq -(\tan(\beta)\,x,\  g5:\ y\sspeq -\tan(\beta\sspm \alpha)\, x }.\pbn
\pbn
\psn
\pbn\pbn

\pbn
\pbn
Keywords:  Origami mathematics, {\sl  Haga}'s theorems, cubic equation, heptagon equation, doubling the cube, trisecting an angle.\psn
AMS MSC numbers: 94B27, 11D25, 97G70 \psn
OEIS A-numbers: \seqnum{A002580},  \seqnum{A246644}.
\end{document}